\documentclass{amsart}
\usepackage{multirow}
\usepackage{booktabs,caption}
\usepackage{mathtools}
\usepackage{amsthm}
\usepackage{amsmath}
\usepackage{graphicx}
\usepackage{color}
\usepackage{amsfonts}
\usepackage{amssymb}
\usepackage{mathrsfs}
\usepackage{amscd}
\usepackage{url}
\usepackage{xcolor}
\usepackage{diagbox}
\usepackage{makecell}
\usepackage[numbers,sort&compress]{natbib}
\usepackage{enumitem}
\usepackage{multirow}
\usepackage{pdflscape}
\usepackage{rotating}

\newtheorem{theorem}{Theorem}[section]
\newtheorem{thm}[theorem]{Theorem}
\newtheorem{prop}[theorem]{Proposition}
\newtheorem{lem}[theorem]{Lemma}

\newtheorem{alg}[theorem]{Algorithm}

\newtheorem{Example}[theorem]{Example}

\newtheorem{rmk}[theorem]{Remark}

\newtheorem{cor}[theorem]{Corollary}

\numberwithin{equation}{section}

\numberwithin{equation}{section}

\newcommand{\re}{\mathbb{R}}

\newcommand{\N}{\mathbb{N}}

\newcommand{\bdS}{\mathcal{S}}

\newcommand{\lmd}{\lambda}

\newcommand{\eps}{\epsilon}

\def\af{\alpha}
\def\bt{\beta}
\def\gm{\gamma}
\def\rank{\mbox{rank}}

\definecolor{green}{rgb}{0.00,0.50,0.00}

\newcommand{\sig}{\sigma}

\newcommand{\reff}[1]{(\ref{#1})}

\newcommand{\mc}[1]{\mathcal{#1}}

\newcommand{\st}{\mathit{s.t.}}
\newcommand{\Span}{\mathrm{Span}}

\newcommand{\mt}[1]{\mathtt{#1}}
\newcommand{\tr}[1]{\mathrm{Trace}(#1)}

\newcommand{\ve}{\Vert}

\newcommand{\bdes}{\begin{description}}
\newcommand{\edes}{\end{description}}

\newcommand{\bal}{\begin{align}}
\newcommand{\eal}{\end{align}}

\newcommand{\bnum}{\begin{enumerate}}
\newcommand{\enum}{\end{enumerate}}

\newcommand{\bit}{\begin{itemize}}
\newcommand{\eit}{\end{itemize}}

\newcommand{\bea}{\begin{eqnarray}}
\newcommand{\eea}{\end{eqnarray}}
\newcommand{\be}{\begin{equation}}
\newcommand{\ee}{\end{equation}}

\newcommand{\baray}{\begin{array}}
\newcommand{\earay}{\end{array}}

\newcommand{\bsry}{\begin{subarray}}
\newcommand{\esry}{\end{subarray}}

\newcommand{\bca}{\begin{cases}}
\newcommand{\eca}{\end{cases}}

\newcommand{\bcen}{\begin{center}}
\newcommand{\ecen}{\end{center}}

\newcommand{\bbm}{\begin{bmatrix}}
\newcommand{\ebm}{\end{bmatrix}}

\newcommand{\bmx}{\begin{matrix}}
\newcommand{\emx}{\end{matrix}}

\newcommand{\bpm}{\begin{pmatrix}}
\newcommand{\epm}{\end{pmatrix}}

\newcommand{\btab}{\begin{tabular}}
\newcommand{\etab}{\end{tabular}}

\begin{document}

\title[A Tight SDP Relaxation for the CQR problem]
{A Tight SDP Relaxation for the Cubic-Quartic Regularization Problem}

\author[J.~Zhou]{Jinling Zhou}
\address{Jinling~Zhou,
School of Mathematics and Computational Science,
Xiangtan University, Xiangtan, Hunan, 411105, China.}
\email{jinlingzhou@smail.xtu.edu.cn}

\author[X.~Liu]{Xin Liu}
\address{Xin Liu,
Academy of Mathematics and Systems Science,
Chinese Academy of Sciences, Beijing 100190, China,
and School of Mathematical Sciences,
University of Chinese Academy of Sciences, Beijing 100049, China.}
\email{liuxin@lsec.cc.ac.cn}

\author[J.~Nie]{Jiawang Nie}
\address{Jiawang Nie,  Department of Mathematics,
University of California San Diego,
9500 Gilman Drive, La Jolla, CA, USA, 92093.}
\email{njw@math.ucsd.edu}

\author[X.~Tang]{Xindong Tang}
\address{Xindong~Tang, Department of Mathematics,
Hong Kong Baptist University,
Kowloon Tong, Kowloon, Hong Kong.}
\email{xdtang@hkbu.edu.hk}

\begin{abstract}
This paper studies how to compute global minimizers of the
cubic-quartic regularization (CQR) problem
\[
\min_{s \in \mathbb{R}^n}  \quad
f_0+g^Ts+\frac{1}{2}s^THs+\frac{\beta}{6}\| s \|^3+
\frac{\sigma}{4} \| s\|^4,
\]
where $f_0$ is a constant, $g$ is an $n$-dimensional vector,
$H$ is an $n$-by-$n$ symmetric matrix, and $\| s \|$ denotes the Euclidean norm of $s$.
The parameter $\sigma$ is nonnegative while $\beta$ can have any sign.
The CQR problem arises as a critical subproblem for
getting efficient regularization methods
for solving unconstrained nonlinear optimization.
Its properties are recently well studied by Cartis and Zhu
{\it [cubic-quartic regularization models for solving
polynomial subproblems in third-order tensor methods, Math. Program, 2025]}.
We propose a structured semidefinite programming (SDP) relaxation method
for solving the CQR problem globally.
The SDP relaxation has only three symmetric positive semidefinite matrix variables
of sizes $(n+1)$-by-$(n+1)$, $3$-by-$3$ and $2$-by-$2$ respectively.
We show that our SDP relaxation is tight if and only if
$\| s^* \| ( \beta + 3 \sigma \| s^* \|) \ge 0$ holds for a global minimizer $s^*$.
When $s^* \ne 0$, this aligns with the sufficient global optimality condition
$\beta + 3 \sigma \| s^* \| \ge 0$ given by Cartis and Zhu.
In particular, if either $\beta \ge 0$ or $H$ has a nonpositive eigenvalue,
then the SDP relaxation is shown to be tight.
Second, we show that all nonzero global minimizers have the same Euclidean norm
for the tight case.
Third, we give an algorithm to detect tightness
and to obtain the set of all global minimizers.
Numerical experiments demonstrate that our SDP relaxation method
is both effective and computationally efficient.
This paper gives a polynomial time algorithm for solving the CQR problem globally,
under the sufficient global optimality condition.
\end{abstract}



\maketitle

\section{Introduction}

In this paper, we focus on computing global minimizers of
the following unconstrained optimization problem
\be  \label{CQR:POLY}
\min\limits_{s \in \re^n} \quad M(s)  \, \coloneqq \,
f_0+g^Ts+\frac{1}{2}s^THs+\frac{\beta}{6}\ve s \ve ^3+\frac{\sigma}{4} \ve s\ve^4,
\ee
where the symmetric matrix $H\in \re^{n\times n}$, the vector
$g\in \re^{n}$, the constant $f_0\in \re$ are problem parameters,
the constants $\beta \in \re$ and $\sigma\ge 0$ are regularization parameters,
the vector $s \in \re^n$ stands for the decision variables and
its Euclidean norm is denoted as $\| s \| = \sqrt{s^Ts}$.
Throughout the paper, we assume $\beta >0$ for the case $\sigma = 0$,
to ensure the existence of global minimizers.
The problem~\reff{CQR:POLY} is called the cubic-quartic regularization (CQR)
problem \cite{zhu2025CQR}. When the regularization parameter $\sigma$ equals zero,
the problem~\reff{CQR:POLY} is reduced to
a cubic regularization problem \cite{Nesterov2006cubic}.
For a minimizer of \reff{CQR:POLY}, we mean it is a global minimizer
(unless its meaning is otherwise specified).

The cubic-quartic regularization problem is proposed and well studied
by Cartis and Zhu \cite{zhu2025CQR}.
In fact, they consider the following more general case:
\be \label{subproblem:mc}
\min\limits_{s \in \re^n}  \quad
 f_0+g^T s+\frac{1}{2}s^THs+\frac{\beta}{6}\ve s\ve ^3_W+\frac{\sigma}{4}\ve s\ve ^4_W,
\ee
where the norm $\ve s \ve_W = \sqrt{s^T W s }$ is defined by a given
symmetric positive-definite matrix $W$. When $W = I_n$ (the $n$-by-$n$ identity matrix),
the problem \reff{subproblem:mc} is reduced to \reff{CQR:POLY}.
If in addition $\beta=0$, the problem \reff{subproblem:mc}
turns into a quadratic-quartic regularization (QQR) problem \cite{zhu2023QQR}.
After imposing a linear transformation $s=W^{-\frac{1}{2}}\tilde{s}$,
the problem \reff{subproblem:mc} turns into the standard form \reff{CQR:POLY},
resulting from the fact that $\ve s\ve _{W} = \ve \tilde{s}\ve$.
Therefore, it suffices to consider the CQR problem in the form of \reff{CQR:POLY}.

\subsection{Motivations}

The Newton trust-region (NTR) methods are a class of practically useful high-order approaches
for computing second order stationary points of the unconstrained optimization problem
\be \label{minf(x)}
\min_{x\in\re^{n}}  f(x),
\ee
where $f\,:\re^n\to \re$ is a nonlinear smooth function. When the Hessian of $f$ is available,
the global convergence and local quadratic convergence rate can be well established under mild conditions. However, its worst case evaluation complexity is hard to establish. To this end, the so-called adaptive
regularization methods are proposed.
Nesterov and Polyak \cite{Nesterov2006cubic} proposed the cubic
regularization methods (AR2/ARC) and established the worst-case evaluation complexity estimates. In each iteration,
it requires to solve the following cubic regularization subproblem
\be   \label{cubic-newton}
\min_{s\in \re^n}  \quad
f_0+g^Ts+\frac{1}{2}s^THs+\frac{\beta}{6}\ve s \ve ^3,
\ee
for a parameter $\beta>0$. Clearly, \eqref{cubic-newton}
is a special case of \eqref{CQR:POLY} with
$\sigma=0$. Later, Nesterov \cite{Nesterov2008accelearting} proposed
an acceleration technique to ARC for solving convex optimization problems
and achieved better worst-case evaluation complexity
when the subproblem \eqref{cubic-newton} is globally solved.
Moreover, Cartis, Gould and Toint
\cite{cartis2010complexity,cartis2011adaptive1,cartis2011adaptive2}
developed the adaptive cubic regularization (ARC) framework,
proving global convergence and worst-case evaluation complexity results for ARC and its variants.
Their analysis clarified how to adaptively select the regularization parameter $\beta$.
Cartis, Gould and Toint \cite{cartis2022evaluation}
gave a characterization for the global minimizer of the ARC subproblem.
Various numerical methods have been developed to solve it.
We refer to \cite{Lieder2020,Jia2022Solving}
for generalized eigenvalue methods and \cite{Gould2012,Gould2020,cartis2011adaptive1,cartis2011adaptive2}
for Krylov subspace methods.

As a standard cubic regularization problem,
ARC has motivated further development of general AR$p$ methods \cite{Birgin2017worst,cartis2020concise,cartis2020sharp,Cartis2024Efficient}.
For problem (\ref{minf(x)}), let $p\ge 2$ be a power
and assume the $p$th order derivatives of $f$ are available.
The AR$p$ scheme requires to minimize the following function in $s$
at the $k$th iteration point $x^{(k)}$:
\[
m_{p}(s) \coloneqq T_p(x^{(k)},s) + \frac{\sigma_k}{p+1}\Vert s\Vert^{p+1},
\]
where $T_p(x^{(k)})$ is the $p$th order Taylor expansion of $f$ at $x^{(k)}$,
and $\sigma_k > 0$ is the regularization parameter.
When $p = 2$, the AR$p$ method reduces to the ARC method.
A key advantage of the AR$p$ method is that using higher-order derivatives
can improve the worst-case evaluation complexity.
Under the Lipschitz continuity assumption on the $p$th order derivatives of $f$,
the AR$p$ method can find a point satisfying the first order optimality conditions
(up to the tolerance $\eps_{1}$) and satisfying
the second order optimality conditions (up to the tolerance $\eps_{2}$),
after at most $O\big(\max\big(\eps_{1}^{-(p+1)/p},\eps_{2}^{-(p+1)/(p-1)} \big) \big)$ evaluations of $f$ and its derivatives \cite{cartis2020concise,cartis2022evaluation}.
Despite this theoretical advantage,
the AR$p$ subproblems are typically difficult to solve.
Indeed, the AR$p$ subproblems of minimizing $m_p(s)$ are NP-hard even for $p = 3$ (see \cite{luo2010quadratic,burachik2021steklov,zhu2025CQR}).
We refer to \cite{Nesterov2021inexact,Nesterov2023inexact,
Nesterov2025quartic,Nesterov2021implementable,Nesterov2021superfast}
for a family of second-order methods that solve convex AR3 problems.
Moreover, it is proposed in \cite{Ahmadi2024Higher} a higher order regularization method
where the regularization term is chosen such that the subproblem is SOS-convex.
Recently, the quadratic-quartic regularization (QQR) method
is proposed in \cite{zhu2023QQR,Zhu2022Quartic}
for minimizing quartically-regularized cubic polynomials.
For AR3 problems, the QQR method approximates the cubic tensor term
with a linear combination of a quadratic form and a quartic term,
making the inner subproblems globally solvable in general.
Interestingly, the recent work \cite{Zhu2026Sufficiently} shows that
the SOS relaxation is tight for solving AR$3$ problems
when the regularization parameters are sufficiently large.
Also, we remark that AR$p$ is a polynomial optimization problem when $p$ is odd,
which can be efficiently solved
by Moment-SOS methods when the problem size is not large.
We refer to \cite{HKL20,LasBk15,mpo,Lau09}
for introductions to polynomial optimization.

The cubic-quartic regularization (CQR) method \cite{zhu2025CQR} is a further new technique for solving AR3 problems and general unconstrained optimization problems like (\ref{minf(x)}). In the scheme of CQR methods, the parameter $\beta$ for the cubic term in (\ref{CQR:POLY}) provides a scalar approximation for the third-order derivatives at the current iterate, while the parameter $\sigma$ is adjusted to control the regularization and progress of the algorithm. Notably, for AR3 problems, the scheme of CQR methods is guaranteed to find a point satisfying the first-order optimality condition within a tolerance $\eps$ after $O(\eps^{-3/2})$ evaluations for the objective function and its derivatives \cite[Theorem~3.4]{zhu2025CQR}.
Moreover, this complexity bound can be further improved
when the AR3 problem has some special structures.
Furthermore, global optimality conditions for
the nonconvex CQR problem (\ref{CQR:POLY}) are studied in \cite{zhu2025optimal}.
Based on these conditions, a Newton-Cholesky method is
proposed to find minimizers for CQR problems, by solving the secular equation; see \cite[Algorithm~1]{zhu2025optimal}.
Other existing methods could also be adapted for solving CQR subproblems,
including generalized eigenvalue methods \cite{Lieder2020,Jia2022Solving,Adachi2017Solving},
Krylov subspace methods \cite{Gould2012,Gould2020,cartis2011adaptive1,cartis2011adaptive2}
and tensor methods for general nonlinear programs \cite{Chen2026Tensor,He2023Homogeneous}.
However, to solve the CQR problem globally,
it still remains a challenge for developing methods
that are both theoretically guaranteed in finding global minimizers
and that are also efficient in computational practice.

In the literature, there exists much work for solving polynomial optimization problems,
whose objective and constraining functions are given by polynomials.
The Moment-SOS hierarchy is efficient for solving them.
We refer to \cite{HKL20,LasBk15,mpo,Lau09} for introductions to this area.
There exist tight relaxation methods for solving general polynomial optimization problems. For instance, tight Moment-SOS relaxations can be obtained with
gradient ideals \cite{NDS06}, optimality conditions \cite{nieopcd},
and Lagrange multipliers \cite{Tight18}.
We also refer to \cite{mpo,Slot2025} for convex polynomial optimization.
However, these tight relaxations may require computation with polynomials of high degrees. The computational cost grows quickly as the degree increases.

The CQR problem is not a standard polynomial optimization problem,
due to the regularization term.
But it can be reformulated as constrained polynomial optimization
(see Remark~\ref{rmk:pop_reform}).
However, doing this is not computationally preferable,
since it becomes a polynomial optimization problem of degree four.
It is NP-hard to globally optimize
quadratic polynomials under quadratic constraints
\cite{Ahmadi2024Higher,Ahmadi2022Complexity,Murty1987NP}.
Moreover, the problem of globally minimizing cubic forms over unit spheres
is also NP-hard \cite{Nesterov2003}.
The computational complexity of Moment-SOS relaxations grows quickly
as the number of variables increases, which makes it difficult to solve large-scale problems. Recently, a global convergence rate for SOS type approximations
with nonconvex adaptive regularization
was established by Cartis and Zhu~\cite{zhu2024sos}.
It introduces an algorithmic framework that combines the SOS Taylor
model with adaptive regularization techniques
for solving nonconvex smooth optimization problems.

It is an interesting question to efficiently compute global optimizers of the CQR problem~\eqref{CQR:POLY}.
Cartis and Zhu \cite{zhu2025CQR} gave necessary and sufficient global optimality conditions.
How to use these conditions to find global minimizers still remains challenging.
To achieve the best convergence properties of CQR methods, it is highly wanted that
global minimizers of CQR problems can be computed efficiently.

\subsection{Contributions}

To address the challenge of computing global minimizers for the CQR problem~\eqref{CQR:POLY},
we propose a structured semidefinite programming (SDP) relaxation method.
For the CQR problem~\eqref{CQR:POLY} with $n$ variables,
our SDP relaxation has three symmetric positive semidefinite (psd) matrix variables
of sizes $(n+1)$-by-$(n+1)$, $3$-by-$3$ and $2$-by-$2$ respectively.
We prove that this relaxation is tight under a very general assumption,
which indeed aligns with the sufficient global optimality condition
given by Cartis and Zhu in \cite{zhu2025CQR}.
The tightness means that it recovers both the global minimum value
and the global minimizers of \eqref{CQR:POLY}.
Under this general assumption, we give a method that guarantees to
find all global optimizers of the CQR problem~\eqref{CQR:POLY} in polynomial time.
We remark that this SDP relaxation method can be applied to solve CQR subproblems
in adaptive regularization methods (e.g., Algorithms~2-5 in \cite{zhu2025CQR}).
The condition for tight SDP relaxations can be used as guidelines
for parameter selections in developing numerical optimization methods.

More specifically, our main contributions are fourfold:

\begin{enumerate}

\item[i)] We propose the following SDP relaxation for solving \eqref{CQR:POLY}:
\be \label{momentSDP}
\left\{
\baray{cl}
\min\limits_{Y, Z_1, Z_2}   &
\displaystyle
  \bbm f_0 & g^T / 2 \\ g/2 & H/2 \ebm
   \bullet  Y +\frac{\beta}{6} (Z_2)_{22}+\frac{\sigma}{4} (Z_1)_{33}  \\
\st  &Y_{00}=1,\quad (Z_1)_{11}=Y_{00}, \\
\quad & (Z_1)_{12}=(Z_2)_{11},\quad (Z_1)_{22}=(Z_2)_{12},\\
\quad & (Z_1)_{13}=(Z_2)_{12},\quad (Z_1)_{23}=(Z_2)_{22},\\
\quad & (Z_1)_{13} = (Z_1)_{22}=Y_{11}+\cdots +Y_{nn},\\
\quad & Y\in \bdS_{+}^{n+1},\quad Z_1\in\bdS_{+}^{3},\quad Z_2\in\bdS_{+}^{2}.
\earay
\right.
\ee
In the above, the notation $\mc{S}_+^m$ denotes the cone of $m$-by-$m$
symmetric positive semidefinite (psd) matrices
and ``$\bullet$" denotes the Euclidean inner product for matrices.
The label indices of the matrix variable $Y$ are $0,1,\ldots,n$,
and those of $Z_1, Z_2$ are $1,2,3$ and $1,2$, respectively.

\item [ii)]
We establish tightness results for the SDP relaxation \reff{momentSDP}.
It is said to be {\it tight} if the optimal value $\vartheta^*$ of \reff{momentSDP}
equals the optimal value $\mu^*$ of \reff{CQR:POLY}.
First, we prove that $\vartheta^* = \mu^*$ if and only if
the regularization parameters $\bt, \sigma$ satisfy the inequality
\be \label{cond:bt}
\| s^* \| ( \bt  + 3 \sig \| s^* \| ) \ge 0,
\ee
where $s^*$ is a global minimizer of \reff{CQR:POLY}.
When $s^* \ne 0$, the above implies that the SDP relaxation \reff{momentSDP}
is tight if and only if the sufficient optimality condition
$\bt  + 3 \sig \| s^* \| \ge 0$,
given in \cite{zhu2025CQR} (also see Theorem~\ref{thm:opcd2}),
holds at $s^*$.
Moreover, we also show that the inequality \reff{cond:bt}
must hold and the SDP relaxation is tight if either $\bt \ge 0$
or $H$ has a nonpositive eigenvalue.
Moreover, we prove that all nonzero global minimizers of \reff{CQR:POLY}
must have the same Euclidean norm when the SDP relaxation \reff{momentSDP} is tight.

\item [iii)]
Our approach solves the CQR problem \eqref{CQR:POLY}
via the SDP relaxation \eqref{momentSDP}. When it is tight,
we not only obtain the true global minimum value,
but also obtain all the global minimizers.
This is summarized in Algorithm~\ref{alg:finds},
which is based on the dual problem of \eqref{momentSDP}.
Furthermore, the algorithm can also detect cases
where the relaxation is not tight,
thereby providing a built-in check for tightness.

\item [iv)]
Numerical experiments demonstrate that our approach is both effective and computationally efficient.
This provides a polynomial time algorithm for globally solving CQR problems
under conditions for tightness.
We also give a method for detecting whether the relaxation is tight or not,
and  we can obtain all global minimizers if it is tight.
Note that \reff{momentSDP} is a semidefinite program.
The matrix variable $Y$ is $(n+1)$-by-$(n+1)$, while
$Z_1$ (resp., $Z_2$) is $3$-by-$3$ (resp., $2$-by-$2$).
There are totally $8$ equality constraints.
The SDP relaxation \reff{momentSDP} can be solved
in $O(n^{3.5} \ln(1/\eps))$ arithmetic operations
by interior-point methods \cite{Todd01}.
In our experiments, we apply the software {\tt Mosek}
to solve \reff{momentSDP}. It works very efficiently.
For instance, when the dimension $n=1000$,
\reff{momentSDP} can be solved within around one minute.

\end{enumerate}

The rest of this paper is organized as follows.
In Section~\ref{sc:pre}, we give a brief survey on polynomial optimization and
global optimality conditions for the cubic-quartic regularization problem.
In Section~\ref{sc:SDr}, we show how to construct the SDP relaxations
and investigate their properties.
In Section~\ref{sc:tight}, we prove tightness of the SDP relaxations.
In Section~\ref{sc:alg}, we propose an algorithm for detecting tightness
of the SDP relaxations and computing all global minimizers.
The numerical experiments to illustrate the effectiveness and efficiency of our method
for solving the CQR problem \reff{CQR:POLY} are reported in Section~\ref{sc:num}.
Finally, some conclusions are drawn in Section~\ref{sc:con}.

\section{Preliminaries}
\label{sc:pre}

\subsection*{Notation}
The symbol $\N$ (resp., $\re$) denotes the set of nonnegative integers (resp., real numbers). The superscript $^T$ denotes the transpose of a matrix or a vector.
Let $\mc{S}^N$ denote the space of all $N$-by-$N$ real symmetric matrices.
For a symmetric matrix $M$, $M \succeq  0$ (resp., $M\succ 0$)
means that $M$ is positive semidefinite (resp., positive definite),
while $M \prec 0$ means that $-M$ is positive definite.
A matrix is indefinite if it is neither positive semidefinite nor negative semidefinite.
The cone of all $N$-by-$N$ real symmetric positive semidefinite matrices is denoted as $\mc{S}^N_+$.
For $x\in\re^n$, its Euclidean norm is $\ve x\ve=\sqrt{x^Tx}$.
For two matrices $A,\,B \in \re^{m\times n}$,
their inner product is denoted as $A \bullet B = \tr{AB^T}$.
The $\mathbf{1}_n$ (resp., $\mathbf{0}_n$) stands for
the all-one (resp., all-zero) vector of length $n$.
The all-one matrix $\mathbf{1}_{m \times n}$
and all-zero matrix $\mathbf{0}_{m \times n}$ are similarly defined.

\subsection{Some basics about polynomials}
\label{ssc:basic}

For $x=(x_1,\cdots,x_n)$ and $\alpha=(\alpha_1, \ldots, \alpha_n)\in\mathbb{N}^n$,
we denote the monomial power
$
x^{\alpha} \coloneqq  x_1^{\alpha_1}\cdots x_n^{\alpha_n}.
$
The degree of $\af$ is
$\left|\alpha  \right|  \coloneqq  \alpha_1+\cdots+\alpha_n$.
For a degree $d > 0$, we denote the power set
\[
\mathbb{N}_d^n \coloneqq  \left\{ \alpha \in \N^n:\, \left| \alpha \right| \leq d\right\}.
\]
The column vector of all monomials in $x$ and of degrees up to $d$ is denoted as
\[
\left [  x\right ]_{d}  \, \coloneqq  \,
(x^\af)_{ \af \in \N^n_d }.
\]
The length of the vector $\left [  x\right ]_{d}$ is $(\substack{n+d\\d})$.
The notation $\mathbb{R}[x]$ denotes the ring of polynomials in
$x=(x_1,\cdots,x_n)$ and with real coefficients.
For a degree $d$, $\mathbb{R}\left [  x\right ]_d$
denotes the set of polynomials in $\mathbb{R}[x]$ with degrees at most $d$.

A polynomial $f\in\mathbb{R}[ x ]$ is said to be a sum of squares (SOS)
if there exist polynomials $f_1,\cdots,f_k\in\mathbb{R}\left [  x\right ]$
such that $f=f_1^2+\cdots+f_k^2$.
The set of all SOS polynomials in $\re[x]$ is denoted as $\Sigma[x]$.
Note that $\Sigma[x]$ is a convex cone in $\re[x]$.
It is interesting to remark that each polynomial in $\Sigma[x]$
is nonnegative everywhere, while the reverse may not be true \cite{mpo}.
For an even degree $2d$, we denote the truncation
\[
\Sigma[ x ]_{2d}   \, \coloneqq \,
\Sigma[ x ] \cap  \re[x]_{2d} .
\]
Moreover, it holds $f \in \Sigma[ x ]_{2d}$
if and only if $f = [x]_d^T X [x]_d$ for some symmetric matrix $X \succeq 0$.

Let $\re\big[\ve x\ve \big]$ denote the set of univariate polynomials in the norm $\ve x \ve$, i.e.,
\[
\re\big[\ve x\ve \big] = \Big \{ \sum_{i=0}^d c_i \ve x \ve^i: c_i \in \re, d \in \N  \Big \}.
\]
Similarly, a polynomial $p( \ve x \ve) \in \re \big[\ve x\ve \big]$ is said to be SOS
if there exist $p_1,\cdots, p_s \in \re\big[ \ve x\ve \big]$
such that $p = p_1^2+\cdots+ p_s^2$.
The set of all SOS polynomials in $\re\big[\ve x\ve \big]$
is denoted as $\Sigma\big[\ve x \ve \big]$.
We remark that a univariate polynomial is nonnegative
everywhere if and only if it is SOS \cite[Chap.~3]{mpo}.
For a degree $d$, we denote the truncation
\[
\Sigma \big[ \ve x \ve \big]_{d}   \, \coloneqq  \,
\Sigma \big[ \ve x \ve \big] \cap  \re \big[\ve x\ve \big]_d .
\]
In this paper, the truncations $\Sigma \big[ \ve x \ve \big]_2$ and $\Sigma \big[ \ve x \ve \big]_4$
are frequently used. Note that
$q(\ve x \ve) = q_0 + q_1 \ve x \ve + q_2 \ve x \ve^2 \in \Sigma\big[ \ve x \ve \big]_2$
if and only if
\[
\bbm  q_0 & \frac{q_1}{2}  \\  \frac{q_1}{2}  &  q_2
\ebm  \succeq 0.
\]
Similarly, $p(\ve x \ve) = p_0 + p_1 \ve x \ve + p_2 \ve x \ve^2 +
p_3 \ve x \ve^3 + p_4 \ve x \ve^4 \in \Sigma\big[ \ve x \ve \big]_4$
if and only if there exists $t \in \re$ such that
\[
\bbm
p_0 & \frac{p_1}{2}  & t \\  \frac{p_1}{2}  &  p_2 - 2t &  \frac{p_3}{2}
\\ t  & \frac{p_3}{2} & p_4
\ebm  \succeq 0.
\]
We refer to \cite[Section~2.4]{mpo} for semidefinite representations for SOS polynomials.

\subsection{The cubic-quartic regularization problem}
\label{sec:pre2}

We review some basic results by Cartis and Zhu \cite{zhu2025CQR}
about the cubic-quartic regularization polynomial.
Consider the general CQR polynomial
\begin{equation}\label{subproblem:zhu}
M_c(s)=f_i+g_i^T s+\frac{1}{2}s^TH_is+
\frac{\beta}{6}\ve s\ve ^3_W+\frac{\sigma_c}{4}\ve s\ve ^4_W,
\end{equation}
where $\sigma_c\ge 0$, $\beta\in\re $,
$W$ is a symmetric positive definite matrix,
and $\ve s \ve_W = \sqrt{s^T W s }$.
For convenience of referencing here, we follow the notation in \cite{zhu2025CQR},
where the CQR polynomials are typically given
by the (high order) derivatives of the objective function at each CQR iterate.
The index $i$ stands for the $i$th iteration,
and the index $c$ of $M_c$ and $\sigma_c$ stands for the CQR method,
to distinguish from the AR$p$ methods.
These indices are not used in other sections of this paper.
The gradient and Hessian of $M_c(s)$ are
\begin{eqnarray*}
\nabla M_c(s) &=& g_i + H_i s + \frac{\bt}{2} \|s\|_W (Ws) + \sig_c \|s\|^2_W (Ws),  \\
\nabla^2 M_c(s) &=& H_i + \frac{\bt}{2} \Big( W\|s\|_W  + \frac{(Ws)(Ws)^T}{\|s\|_W}  \Big)
 + \sig_c \Big( \|s\|_W^2 W  + 2 (Ws)(Ws)^T \Big) .
\end{eqnarray*}
The optimality conditions for global minimizers of $M_c(s)$
are given by Cartis and Zhu \cite{zhu2025CQR} as follows.

\begin{thm} \cite[Theorem~2.1]{zhu2025CQR} \label{thm:opcd}
Let $s_c$ be a global minimizer of $M_c(s) $ over $\re^n$ and let
\begin{equation} \label{bs}
B(s_c) := H_i + \frac{\beta}{2} W \|s_c\|_W + \sigma_c W \|s_c\|_W^2.
\end{equation}
Then, $ s_c$ satisfies
\[
B(s_c) s_c = -g_i, \quad  B(s_c) \succeq 0.
\]
\end{thm}

\begin{thm} \cite[Theorem 2.2]{zhu2025CQR} \label{thm:opcd2}
Let $ B(s_c)$ be defined as in \reff{bs}.
Then, a point $s_c$ is a global minimizer of $M_c(s)$ over $\mathbb{R}^n$
if the following three conditions hold:
\bnum
\item [(1)] $g_i$ is in the range of $B(s_c)$, such that
\[
B(s_c) s_c = \left( H_i + \frac{\beta \|s_c\|_W}{2}  W +
\sigma_c \|s_c\|_W^2 W \right) s_c = -g_i.
\]

\item [(2)] $B(s_c) $ is positive semidefinite:
\[
B(s_c) \succeq 0.
\]

\item[(3)] $\beta $ satisfies
\[
\beta \geq -3 \sigma_c \|s_c\|_W, \quad \text{or equivalently} \quad
\|s_c\|_W \geq \frac{-\beta}{3 \sigma_c}.
\]
\enum
\end{thm}

These conditions are very useful for studying the cubic-quartic regularization problem.
We refer to \cite{zhu2025CQR} for more details.

\section{The semidefinite relaxation}
\label{sc:SDr}

In this section, we give an SDP relaxation for solving the
cubic-quartic regularization problem \eqref{CQR:POLY}.

Recall the cones $\Sigma[s]_2$, $\Sigma\big[ \ve s \ve \big]_4$, $\Sigma\big[ \ve s \ve \big]_2$
introduced in Subsection~\ref{ssc:basic}.
Let $M(s)$ be the objective function in \reff{CQR:POLY} and let $\mu^*$
denote the minimum value of \reff{CQR:POLY}.
For a scalar $\gamma \in \re$, it holds
\[
M(s) -\gamma \in \Sigma \left[ s \right ]_2+
\Sigma\big[ \ve s \ve \big]_4 + \ve s \ve  \Sigma\big[ \ve s \ve \big]_2
\]
if and only if there exist symmetric matrices $X_0, X_1, X_2$ such that
\be \label{sos:eq}
M(s) -\gamma \, = \,  \bbm 1 \\ s \ebm^T X_0 \bbm 1 \\ s \ebm +
\bbm 1 \\ \| s \| \\ \| s \|^2 \ebm^T X_1 \bbm 1 \\ \| s \| \\ \| s \|^2 \ebm  +
\|s \| \bbm 1 \\ \| s \|  \ebm^T X_2 \bbm 1 \\ \| s \|   \ebm ,
\ee
\be \label{sos:psd}
0 \preceq X_0 \in \mc{S}^{n+1}, \quad
0 \preceq X_1 \in \mc{S}^{3}, \quad
0 \preceq X_2 \in \mc{S}^{2}.
\ee
Note that \reff{sos:eq} implies $M(s) \ge \gamma$ for all $s \in \re^n$,
so $\gamma \le \mu^{*}$.

We consider the following relaxation for solving \reff{CQR:POLY}:
\be\label{problem:dual:sos}
\left\{
\baray{cl}
\max & \quad \gamma \\
\st  & \quad  \gamma \,\, \text{satisfies} \,\, \reff{sos:eq}-\reff{sos:psd} .
\earay
\right.
\ee
Let $\gamma^{*}$ denote the optimal value of \reff{problem:dual:sos}, then
$
\gamma^*  \le   \mu^* .
$
The relaxation \reff{problem:dual:sos} is said to be {\it tight}
if $\gamma^* =  \mu^*$.
The above problem is strictly feasible when $\sigma > 0$.
To see this, let
\[
\begin{gathered}
c_1 \coloneqq \frac{(\beta-6)^2}{36\sigma} + 1,\, c_2 \coloneqq  1-\frac{\lmd_{min}}{2}, \,
c_3 \coloneqq 1+ g^T Q^{-1} g, \, c_4 \coloneqq 1+ w^T P^{-1}w,
\end{gathered}
\]
where
\[ P \coloneqq \begin{bmatrix}  2c_1 & \frac{\beta}{6}-1 \\  \frac{\beta}{6}-1 & \frac{\sigma}{2} \end{bmatrix},\quad
Q \coloneqq H + 2c_2I_n,\quad w \coloneqq \begin{bmatrix}  -1 \\  -c_2 - c_1 \end{bmatrix}.\]
By Schur complement, if we let
\[
X_0 \coloneqq \frac{1}{2} \begin{bmatrix} c_3 & g^T \\ g & Q \end{bmatrix}, \quad
X_1 \coloneqq \frac{1}{2}\begin{bmatrix}  c_4 & w^T \\ w & P \end{bmatrix}, \quad
X_2 \coloneqq I_2, \quad \gm \coloneqq f_0 - \frac{c_3+c_4}{2},
\]
then $X_i\succ 0$ for each $i = 0,1,2$, and $(\gm,X_0,X_1,X_2)$ is feasible for (\ref{problem:dual:sos}).
This implies that (\ref{problem:dual:sos}) is strictly feasible.

In the following, we derive the dual optimization problem of \reff{problem:dual:sos}.
Note $s = (s_1, \ldots, s_n)$ and $M(s)$ belongs to the vector space
\[
V\, \coloneqq  \, \Span \left\{\mathbf{1},s_1,\ldots, s_n, s_1^2,s_1s_2\ldots,s_n^2,
\ve s \ve, \ve s \ve ^2,\ve s \ve^3,\ve s \ve^4 \right\}.
\]
The above $\mathbf{1}$ is the vector of constant one polynomial.
Let $\ell$ be a linear functional on $V$, i.e., a linear mapping from $V$ to $\re$.
Then, $\ell$ is uniquely determined by the values:
\[
\ell( \mathbf{1} ),\ell(s_1),\ldots, \ell(s_n), \ell(s_1^2),
\ell(s_1s_2), \ldots,\ell(s_n^2), \ell(\ve s \ve),\ldots ,\ell(\ve s \ve^4).
\]
Let $Y$, $Z_1$, $Z_2$ be the matrices such that
\begin{eqnarray} \label{mat:Y}
Y &  \coloneqq   & \bbm
\ell( \mathbf{1} )  & \ell(s_1)   &  \ell(s_2)   & \cdots & \ell(s_n)  \\
\ell(s_1)  & \ell(s_1^2) &  \ell(s_1s_2) & \cdots & \ell(s_1s_n)  \\
\ell(s_2)  & \ell(s_1s_2) & \ell(s_2^2) & \cdots & \ell(s_2s_n)  \\
\vdots & \vdots & \vdots  & \ddots & \vdots  \\
\ell(s_n)  & \ell(s_1s_n) & \ell(s_2s_n) & \cdots & \ell(s_n^2)  \\
\ebm \\
& =  & \bbm
Y_{00}  & Y_{01}   &  Y_{02}   & \cdots & Y_{0n}  \\
Y_{10}  & Y_{11}   &  Y_{12}   & \cdots & Y_{1n}  \\
Y_{20}  & Y_{21}   &  Y_{22}   & \cdots & Y_{2n}  \\
\vdots & \vdots & \vdots  & \ddots & \vdots  \\
Y_{n0}  & Y_{n1}   &  Y_{n2}   & \cdots & Y_{nn}  \\
\ebm \in \mc{S}^{n+1},  \nonumber
\end{eqnarray}
\be \label{mat:Z1}
Z_1   \, \coloneqq \,    \bbm
\ell( \mathbf{1} ) & \ell( \ve s \ve )   & \ell( \ve s \ve^2 ) \\
\ell( \ve s \ve)   & \ell( \ve s \ve^2 ) & \ell( \ve s \ve^3 ) \\
\ell( \ve s \ve^2) & \ell( \ve s \ve^3 ) & \ell( \ve s \ve^4 )
\ebm
=\bbm
(Z_1)_{11} & (Z_1)_{12} & (Z_1)_{13} \\
(Z_1)_{21} & (Z_1)_{22} & (Z_1)_{23} \\
(Z_1)_{31} & (Z_1)_{32} & (Z_1)_{33}
\ebm\in \mc{S}^{3},
\ee
\be \label{mat:Z2}
Z_2   \, \coloneqq \,   \bbm
\ell( \ve s \ve ) & \ell( \ve s \ve^2 ) \\
\ell( \ve s \ve^2 ) &  \ell( \ve s \ve^3 )
\ebm
=\bbm
(Z_2)_{11} & (Z_2)_{12} \\
(Z_2)_{21} & (Z_2)_{22}
\ebm \in \mc{S}^{2}.
\ee
The above matrices $Y$, $Z_1$ and $Z_2$ determine
the linear functional $\ell$ on $V$.
For $\ell$ to be well defined on $V$,
the relation $\|s\|^2 = s_1^2 + \cdots + s_n^2$ poses the condition
\[
\ell (\|s\|^2 ) \, = \,  \ell( s_1^2 ) + \cdots + \ell( s_n^2 ).
\]
The value $\ell (\|s\|^2 )$ appears in $(Z_1)_{13}, (Z_1)_{22}, (Z_1)_{31}, (Z_2)_{12}$
and $(Z_2)_{21}$, while $\ell( s_1^2 )$, \ldots, $\ell( s_n^2 )$
appear in the diagonal of $Y$. This implies the equations
\be \label{eq:YZ1Z2}
(Z_1)_{13} =  (Z_1)_{22} =  (Z_1)_{31} = (Z_2)_{12} = (Z_2)_{21} =
Y_{11}  + \cdots + Y_{nn}.
\ee
Indeed, there is a one-to-one correspondence between linear functionals $\ell$ on $V$
and triples $(Y, Z_1, Z_2) \in \mc{S}^{n+1}\times \mc{S}^{3}\times \mc{S}^{2}$
satisfying (\ref{eq:YZ1Z2}).
Denote the cone
\[
K  \, \coloneqq \,   \Sigma [s]_2 +
\Sigma\big[ \ve s \ve \big]_4 + \ve s \ve  \Sigma\big[ \ve s \ve \big]_2
\subseteq V.
\]

\begin{lem}  \label{lem:dual}
Let $\ell$ be a linear functional on $V$ and let
$Y$, $Z_1$, $Z_2$ be the matrices as in
\reff{mat:Y}, \reff{mat:Z1}, \reff{mat:Z2} respectively.
Then,  $\ell \ge 0$ on $K$ if and only if
\[
Y \succeq 0, \quad Z_1 \succeq 0, \quad Z_2 \succeq 0.
\]
\end{lem}
\begin{proof}
For a set $T \subseteq V$, its conic hull is denoted as
\[
\mbox{cone}(T) \, \coloneqq \,  \Big \{ \sum_{i=1}^N \lmd_i t_i :
\lmd_i \ge 0, t_i \in T, N \in \N \Big \}.
\]
Observe that
\[
\Sigma[s]_2  = \mbox{cone}\{p_1^2: \, p_1 \in \re[s]_1 \},  \quad
\Sigma\big[ \ve s \ve \big]_4  =
\mbox{cone}\{p_2^2: \, p_2 \in \re\big[\Vert s\Vert \big]_2 \},
\]
\[
\ve s \ve  \Sigma\big[ \ve s \ve \big]_2 =  \mbox{cone}
\{\ve s \ve p_3^2: \, p_3\in \re\big[\Vert s\Vert\big]_1 \} .
\]
So $\ell \ge 0$ on $K$ if and only if
\[
\ell(p_1^2) \ge 0, \quad  \ell(p_2^2) \ge 0, \quad  \ell(p_3^2 \ve s \ve ) \ge 0,
\]
for all $p_1 \in \re[s]_1$, $p_2 \in \re\big[ \ve s \ve \big]_2$ and
$p_3 \in \re\big[\ve s \ve \big]_1$. One can write that
\[
p_1 = v_1^T \bbm 1 \\ s \ebm, v_1 \in \re^{n+1}, \quad
p_2 = v_2^T \bbm 1 \\ \| s \| \\  \| s \|^2  \ebm, v_2 \in \re^{3}, \quad
p_3 = v_3^T \bbm 1 \\ \| s \|  \ebm, \, v_3 \in \re^{2} .
\]
Observe the relations
\[
\ell(p_1^2) = v_1^T Y v_1, \quad  \ell(p_2^2) = v_2^T Z_1 v_2, \quad
\ell(p_3^2 \ve s \ve ) = v_3^T Z_2 v_3.
\]
Then the conclusion of the lemma follows from the above.
\qed
\end{proof}

Let $\ell$ be a linear functional on $V$, and let $Y,Z_1,Z_2$ be matrices given as in (\ref{mat:Y})-(\ref{mat:Z2}).
Define the linear function
\[
 \vartheta(Y, Z_1, Z_2) \coloneqq  \ell ( M(s) ) =
 \bbm
f_0 & g^T / 2 \\
g/2 & H/2
\ebm
 \bullet  Y +\frac{\beta}{6} (Z_2)_{22}+\frac{\sigma}{4} (Z_1)_{33} .
\]
In the above, the bullet $\bullet$ denotes the Euclidean inner product (see Section~\ref{sc:pre}).
It is interesting to note that
\[
\vartheta(Y, Z_1, Z_2) - \gamma  \, = \, \ell ( M(s) -\gamma \mathbf{1})
-\gamma(1- Y_{00}).
\]
Thus, if $Y_{00} = 1$, $\ell \ge 0$ on $K$,
and $\gamma$ is feasible for \reff{problem:dual:sos},
then $\vartheta(Y, Z_1, Z_2) \ge \gamma$.
By Lemma~\ref{lem:dual}, the dual optimization problem of
the relaxation~\reff{problem:dual:sos} is
\be \label{pro-moment}
\left\{
\baray{cl}
\min\limits_{Y, Z_1, Z_2}   &   \vartheta(Y, Z_1, Z_2)   \\
\st  &Y_{00}=1,\quad (Z_1)_{11}=Y_{00}, \\
\quad & (Z_1)_{12}=(Z_2)_{11},\quad (Z_1)_{22}=(Z_2)_{12},\\
\quad & (Z_1)_{13}=(Z_2)_{12},\quad (Z_1)_{23}=(Z_2)_{22},\\
\quad & (Z_1)_{13} = (Z_1)_{22}=Y_{11}+\cdots +Y_{nn},\\
\quad & Y\in \bdS_{+}^{n+1},\quad Z_1\in\bdS_{+}^{3},\quad Z_2\in\bdS_{+}^{2}.
\earay
\right.
\ee
Note that \reff{pro-moment} is an SDP problem.
The matrix variable $Y$ is the dual variable
for the psd constraint $X_0\succeq 0$ in (\ref{problem:dual:sos}),
while $Z_1$ and $Z_2$ are dual variables for constraints
$X_1 \succeq 0$ and $X_2 \succeq 0$ respectively.
One may check that $Y$ is $(n+1)$-by-$(n+1)$, and
$Z_1$ (resp., $Z_2$) is $3$-by-$3$ (resp., $2$-by-$2$).
There are totally $8$ equality constraints.
Therefore, the relaxation \reff{momentSDP} can be solved
in $O(n^{3.5} \ln(1/\eps))$ arithmetic operations
by path-following interior-point methods \cite{Todd01}.

In the following, we provide another perspective
of the relaxation (\ref{pro-moment}).
Let $Y\in \mc{S}^{n+1}$, $Z_1\in \mc{S}^{3}$ and $Z_2\in \mc{S}^{2}$ be matrices.
Note that there exists a point $s\in \re^n$ such that
\be\label{eq:YZ1Z2=ssT}
\begin{gathered}
Y = \bbm 1  & s^T \\ s   &  ss^T \ebm , \quad
Z_1 = \bbm
1             &   \ve s \ve    &   \ve s \ve^2  \\
\ve s \ve   &   \ve s \ve^2  &   \ve s \ve^3  \\
\ve s \ve^2 &   \ve s \ve^3  &   \ve s \ve^4  \\
\ebm , \quad
Z_2 = \bbm
\ve s \ve   &   \ve s \ve^2    \\
\ve s \ve^2 &   \ve s \ve^3    \\
\ebm
\end{gathered}
\ee
if and only if $Y_{00} = 1$, (\ref{eq:YZ1Z2}) holds, and
\[
Y \succeq 0,\quad Z_1\succeq 0,\quad Z_2\succeq 0,\quad
\rank\, Y = \rank\, Z_1 = \rank\, Z_2 = 1.
\]
Clearly, if (\ref{eq:YZ1Z2=ssT}) holds, then
$
M(s) = \vartheta(Y, Z_1, Z_2 ).
$
Thus, the CQR problem (\ref{CQR:POLY})
can be equivalently reformulated as
\be \label{pro-moment-rank1}
\left\{
\baray{cl}
\min\limits_{Y, Z_1, Z_2}   &   \vartheta(Y, Z_1, Z_2)   \\
\st  & Y_{00}=1,\quad (Z_1)_{11}=Y_{00}, \\
\quad & (Z_1)_{12}=(Z_2)_{11},\quad (Z_1)_{22}=(Z_2)_{12},\\
\quad & (Z_1)_{13}=(Z_2)_{12},\quad (Z_1)_{23}=(Z_2)_{22},\\
\quad & (Z_1)_{13} = (Z_1)_{22}=Y_{11}+\cdots +Y_{nn},\\
\quad & \rank\, Y = \rank\, Z_1 = \rank\, Z_2 = 1, \\
\quad & Y\in \bdS_{+}^{n+1},\quad Z_1\in\bdS_{+}^{3},\quad Z_2\in\bdS_{+}^{2}.
\earay
\right.
\ee
The above problem is difficult to solve, due to the rank one constraints.
We can get the relaxation (\ref{pro-moment})
by dropping the rank one constraints in (\ref{pro-moment-rank1}).
In other words, matrix variables $Y$, $Z_1$ and $Z_2$ in (\ref{pro-moment})
can be viewed as relaxations of the matrices in (\ref{eq:YZ1Z2=ssT}).
Clearly, all feasible points of (\ref{pro-moment-rank1})
are feasible for (\ref{pro-moment}), while the reverse may not be true.
So, the optimal value $\vartheta^*$ of \reff{pro-moment}
is less than or equal to the minimum value $\mu^{*}$ of \reff{CQR:POLY}, i.e.,
$\vartheta^*  \le   \mu^*.$
The relaxation \reff{pro-moment} is said to be {\it tight} if $\vartheta^* =  \mu^*$.

The following is a basic property about the relaxations
\reff{problem:dual:sos} and \reff{pro-moment}.
In particular, the relaxation \reff{problem:dual:sos}
is tight if and only if \reff{pro-moment} is tight.

\begin{thm}  \label{thm:bound}
Let $\mu^*$, $\gamma^*$ and $\vartheta^*$ denote the optimal values of
\reff{CQR:POLY}, \reff{problem:dual:sos} and \reff{pro-moment} respectively.
Then, it holds that
\be \label{bd:rel}
\gamma^*  \, = \,  \vartheta^*  \, \le \,  \mu^*.
\ee
Moreover, the optimal value $\gamma^*$ is attainable for \reff{problem:dual:sos}.
\end{thm}
\begin{proof}
Suppose $\gamma$ and $(Y, Z_1, Z_2)$ are feasible points of
\reff{problem:dual:sos} and \reff{pro-moment} respectively.
Since \reff{problem:dual:sos} and \reff{pro-moment} are dual to each other,
we have $\gamma^*  \le  \vartheta^*$ by weak duality.
We show that \reff{pro-moment} has a strictly feasible point.
For $s\in \re^n$, let $Y(s)$, $Z_1(s)$, $Z_2(s)$
be the matrices as in (\ref{eq:YZ1Z2=ssT}).
Let $\nu$ denote the normal distribution on $\re^n$ and
\[
\hat{Y} = \int Y(s) \mt{d} \nu, \quad
\hat{Z}_1 = \int Z_1(s) \mt{d} \nu, \quad
\hat{Z}_2 = \int Z_2(s) \mt{d} \nu.
\]
Since $Y(s)\in \bdS_{+}^{n+1}$, $Z_1(s)\in\bdS_{+}^{3}$,
$Z_2(s)\in\bdS_{+}^{2}$ for all $s\in\re^n$,
$\hat{Y},\hat{Z}_1,\hat{Z}_2$ are all positive definite matrices.
Then $(\hat{Y},\hat{Z}_1,\hat{Z}_2)$
is a strictly feasible point for \reff{pro-moment}.
By the strong duality theorem, we must have
$\gamma^* =  \vartheta^*$ and $\gamma^*$ is  attainable.
The inequality $\vartheta^* \le \mu^{*}$ is shown earlier.
\qed
\end{proof}

\begin{rmk}\label{rmk:pop_reform}
By introducing the auxiliary variable $t$, the CQR problem (\ref{CQR:POLY})
can be equivalently reformulated as the polynomial optimization problem:
\be\label{eq:trans_pop}
\left\{ \begin{array}{cl}
\displaystyle\min_{s \in \re^n, t \in \re}  & f_0+g^Ts+\frac{1}{2}s^THs +
\frac{\beta}{6}t ^3 + \frac{\sigma}{4} t^4, \\
\st & t^2= s^T s,\ t\ge 0.
\end{array}
 \right.
\ee
The Moment-SOS relaxations are
applicable for globally solving \reff{eq:trans_pop}.
Since it has degree four, the order $k$ of Moment-SOS relaxations is at least $2$.
For the $k$th order relaxation, there is one psd matrix variable of
size $\binom{n+1+k}{n+1} \times \binom{n+1+k}{n+1}$
and one psd matrix variable of
size $\binom{n+k}{n+1} \times \binom{n+k}{n+1}$.
We refer to \cite[Chap.~5]{mpo} for the details.
In contrast, in our SDP relaxation (\ref{pro-moment}),
the psd variable $Y$ is $(n+1)$-by-$(n+1)$, and
$Z_1$ (resp., $Z_2$) is $3$-by-$3$ (resp., $2$-by-$2$),
leading to a significantly more scalable method.
\end{rmk}

The SDP relaxation \reff{pro-moment} has totally $8$ equality constraints.
It can be solved in $O(n^{3.5} \ln(1/\eps))$ arithmetic operations
by interior-point methods \cite{Todd01}.
So, the relaxations (\ref{problem:dual:sos}) and (\ref{pro-moment})
are polynomial time solvable.

\begin{prop}  \label{prop:rankone}
Suppose $(Y^*, Z_1^*, Z_2^*)$ is an optimizer of \reff{pro-moment}.
If all $Y^*$, $Z_1^*$, $Z_2^*$ are rank one, then the SDP relaxations
\reff{problem:dual:sos} and \reff{pro-moment} are tight
(i.e., $\gamma^*  =  \vartheta^*  =  \mu^{*}$),
and $s^* \coloneqq (Y^*_{10},\ldots,Y^*_{n0})$
is a global minimizer of $M(s)$.
\end{prop}
\begin{proof}
Suppose $(Y^*, Z_1^*, Z_2^*)$ is an optimizer of \reff{pro-moment}.
Since $Y^*$ is rank one, it holds
\[  Y^* = \bbm 1   \\ s^*  \ebm \bbm 1 & \,\,  {s^*}^T \ebm  \succeq 0 . \]
By the constraints of (\ref{pro-moment}), we also have
\[ (Z_1)_{22} =  (Z_1)_{13} = (Z_2)_{12} = {s^*}^T {s^*} = \Vert s^*\Vert^2. \]
Because both $Z_1^*$ and $Z_2^*$ are rank one and $(Z_1)_{12} = (Z_2)_{11} \ge 0$,
we further have
\be\label{eq:rank-1}\begin{gathered}
Z_1^* = \bbm
1             &   \ve s^* \ve    &   \ve s^* \ve^2  \\
\ve s^* \ve   &   \ve s^* \ve^2  &   \ve s^* \ve^3  \\
\ve s^* \ve^2 &   \ve s^* \ve^3  &   \ve s^* \ve^4  \\
\ebm \succeq 0, \quad
Z_2^* = \bbm
\ve s^* \ve   &   \ve s^* \ve^2    \\
\ve s^* \ve^2 &   \ve s^* \ve^3    \\
\ebm \succeq 0.
\end{gathered}
\ee
Thus, we get
\[
M(s^*) \, = \, \vartheta(Y^*, Z_1^*, Z_2^*)  \, =  \,  \vartheta^*.
\]
Since $\mu^*  \le  M(s^*)$, the above together with \reff{bd:rel} implies
$\mu^*  =  M(s^*)$. This means that
$s^*$ is a global minimizer of $M(s)$,
which completes the proof.
\qed
\end{proof}

Proposition~\ref{prop:rankone} implies that if
$\rank\, Y^* = \rank\, Z_1^* = \rank\, Z_2^* = 1$,
then the SDP relaxation \reff{pro-moment} must be tight,
and we can get a minimizer $s^*$ for the CQR problem (\ref{CQR:POLY}).
However, this does not necessarily mean that this $s^*$
is the unique minimizer of (\ref{CQR:POLY}).
Indeed, when the relaxation (\ref{pro-moment}) is tight,
there always exists an optimizer $(Y^*, Z_1^*, Z_2^*)$ of \reff{pro-moment}
such that these three matrices are rank one, regardless of
uniqueness of the minimizer for the CQR problem.
Also, $Y^*$, $Z_1^*$ and $Z_2^*$ being rank one is
a sufficient but not necessary condition
for the tightness of (\ref{pro-moment}),
i.e., the relaxation may still be tight
when some minimizer matrices are not rank one.
The following is such an exposition.

\begin{Example}\rm
\label{ex:rank1_inf_global_min}
Consider the CQR problem
\be\label{eq:s=0ors=2}
\min_{s\in \re^n} \quad  \Vert s\Vert ^4-4\Vert s \Vert ^3+ 4 \| s \|^2
 = \| s \|^2 (\|s\| - 2)^2.
\ee
The minimum value $\mu^* = 0$. The global minimizers are $\mathbf{0}_n$
and all points $s^*$ such that $\| s^* \| = 2$.
For each $n = 1,2,\ldots$, the SDP relaxation (\ref{pro-moment})
has a minimizer $(Y^*, Z_1^*, Z_2^*)$ such that
$\rank\, Y^* = \rank\, Z_1^* = \rank\, Z_2^* = 1$.
For instance, when $n = 2$, let
\[
Y^* =
\bbm
1  &  \sqrt{2} & \sqrt{2}  \\
\sqrt{2}  &  2 & 2  \\
\sqrt{2}  &  2 & 2  \\
\ebm, \quad
Z_1^*=
\bbm
1   &  2   &  4 \\
2   &  4   &  8  \\
4   &  8   &  16  \\
\ebm, \quad
Z_2^*=
\bbm
2 &  4  \\
4 &  8 \\
\ebm.
\]
Then, one can check that
\[
\vartheta(Y^*, Z_1^*, Z_2^*) \, = \,
(Z_1^*)_{33}-4(Z_2^*)_{22}+4((Y^*)_{22}+(Y^*)_{33}) = 0.
\]
Note that
\be \label{eq:rank1_inf_global_min_sos}
\| s \|^2 (\|s\| - 2)^2 - 0 = (\|s\|^2 - 2\|s\|)^2
\in \Sigma\big[\Vert s\Vert \big]_4 + \Vert s\Vert\Sigma\big[\Vert s\Vert \big]_2.
\ee
By Theorem~\ref{thm:bound}, we have $\gamma^* = \vartheta^* =  0$,
and the above triple $(Y^*, Z_1^*, Z_2^*)$ is an optimizer of (\ref{pro-moment}).
Furthermore, the ranks of $Y^*$, $Z_1^*$ and $Z_2^*$ are all equal to one,
so the relaxation (\ref{pro-moment}) is tight and
$s^* = (\sqrt{2} , \sqrt{2})$ is a minimizer for the CQR problem,
by Proposition~\ref{prop:rankone}.
However, $s^*$ is not the unique minimizer,
since there are infinitely many of them.
Finally, we remark that there also exist minimizer matrices
for (\ref{pro-moment}) that are not rank one.
For instance, if we let
\[
\hat{Y} \, =  \,
\bbm
1  &  \frac{1}{2} & \frac{1}{2}  \\
\frac{1}{2}  &  1 & 0  \\
\frac{1}{2}  &  0 & 1  \\
\ebm, \quad
\hat{Z}_1 \, =  \,
\bbm
1   &  1   &  2 \\
1   &  2   &  4  \\
2   &  4   &  8  \\
\ebm, \quad
\hat{Z}_2 \, = \,
\bbm
1 &  2  \\
2 &  4 \\
\ebm,
\]
then $\vartheta(\hat{Y}, \hat{Z}_1, \hat{Z}_2) = \gamma^* = \mu^* = 0.$
So, $(\hat{Y}, \hat{Z}_1, \hat{Z}_2)$ is also a minimizer of (\ref{pro-moment}) and
$\rank\, \hat{Y} = 3$ and $\rank\, \hat{Z}_1 = 2$,
while the relaxation (\ref{pro-moment}) is tight.
\end{Example}

Interestingly, when the SDP relaxations (\ref{problem:dual:sos}) and (\ref{pro-moment}) are tight, we can get all minimizers for the CQR problem, even if there are infinitely many of them. We will discuss this issue in Section~\ref{sc:alg}.

\section{Tightness of the relaxations}
\label{sc:tight}

In this section, we prove the tightness  of
SDP relaxations \reff{problem:dual:sos} and \reff{pro-moment}
under some general assumptions.
Consider the CQR problem  \reff{CQR:POLY}.
We assume $M(s)$ has a global minimizer $s^*$.
This can be ensured by the assumption that
$\sigma > 0$, or $\sig =0$ but $\beta > 0$,
or $\sig = \beta =0$ but $H$ is positive definite.

In Subsection~\ref{sec:pre2}, when the weight matrix $W$ is identity,
the matrix function in \reff{bs} reduces to
\[
B(s) = H + \frac{\bt}{2}  \|s\| I_n    + \sig   \|s\|^2 I_n .
\]
By Theorem~\ref{thm:opcd} (also see \cite{zhu2025CQR}),
if $s^*$ is a global minimizer of $M(s)$, then it must satisfy
\[
B(s^*) s^* = - g, \quad B(s^*) \succeq 0.
\]
By Theorem \ref{thm:opcd2} (also see \cite{zhu2025CQR}),
the point $s^*$ must be a global minimizer of $M(s)$ if it satisfies
\[
B(s^*) s^* = -g, \quad
B(s^*) \succeq 0, \quad
\bt \ge -3 \sig \| s^* \|.
\]

First, we prove a sufficient and necessary condition for the SDP relaxations
\reff{problem:dual:sos} and \reff{pro-moment} to be tight.
Recall that $\mu^{*}$, $\gamma^*$, $\vartheta^*$ denote the optimal values
of \reff{CQR:POLY}, \reff{problem:dual:sos}, \reff{pro-moment} respectively.
It is worthy to note that \reff{problem:dual:sos} is tight if and only if
\reff{pro-moment} is tight. This is shown in Theorem~\ref{thm:bound}.

\begin{thm}\label{thm:tight_cond}
Let $s^*$ be a global minimizer of \reff{CQR:POLY}.
Then, the SDP relaxations \reff{problem:dual:sos} and \reff{pro-moment}
are tight (i.e., $\gamma^*  =  \vartheta^*  \, = \,  \mu^{*}$)
if and only if
\be \label{bt>=-3}
\| s^* \| \big(  \bt + 3 \sig \| s^* \| \big) \ge 0.
\ee
Furthermore, if \reff{bt>=-3} holds for one minimizer $s^*$,
then it holds for all minimizers.
\end{thm}
\begin{proof}
``$\Leftarrow$"  Assume \reff{bt>=-3} holds,
we need to show that \reff{problem:dual:sos} and \reff{pro-moment}
are both tight.

First, we consider the case that $\|s^*\| = 0$, i.e., $s^* = 0$.
By Theorem~\ref{thm:opcd}, we have $-g = B(s^*)s^* = \mathbf{0}_n$.
Since $\mu^{*} = M(s^*)$, we have
\[
M(s) -\mu^{*} =
\frac{1}{2} s^T H s + \frac{\bt}{6}\Vert s\Vert^3 +
\frac{\sigma}{4} \Vert s\Vert^4 \ge 0\quad \mbox{for all} \quad s\in \re^n.
\]
Let $\lmd_{min}$ denote the smallest eigenvalue of $H$
with the eigenvector $\hat{s}$, then
\[
M(s) -\mu^{*} =
\frac{1}{2} s^T \Big( H - \lmd_{min}I_n \Big) s +
\frac{1}{2} \lmd_{min} \Vert s \Vert^2  +\frac{\bt}{6}\Vert s\Vert^3 +
\frac{\sigma}{4} \Vert s\Vert^4.
\]
In particular, it holds
\[
M(\hat{s}) -\mu^{*} =
\frac{1}{2} \lmd_{min} \Vert \hat{s} \Vert^2  +\frac{\bt}{6}\Vert \hat{s} \Vert^3 +
\frac{\sigma}{4} \Vert \hat{s} \Vert^4.
\]
Note $M(s) -\mu^{*} \ge 0$ for all $s\in \re^n$ and the eigenvector
$\hat{s}$ can be selected such that $\Vert \hat{s}\Vert$ is an arbitrary positive number.
So, we get
\[
\rho(z)  \coloneqq
\frac{1}{2} \lmd_{min} z^2  +\frac{\bt}{6} z^3 +
\frac{\sigma}{4} z^4 \ge 0 \quad \text{for all} \,\, z \ge 0.
\]
By Theorem~3.2.1 of \cite{mpo}, there exist a quadratic polynomial
$p(z)$ and a linear polynomial $q(z)$ such that
$
\rho(z)  =  p( z )^2 + z q(z)^2.
$
Hence,
\[
M(s) - \mu^{*}  = \frac{1}{2} s^T \Big( H - \lmd_{min}I_n \Big) s
+ p( \|s\| )^2 + \|s\| q(\|s\|)^2.
\]
Since $H - \lmd_{min}I_n \succeq 0$, we get
\[
M(s) - \mu^{*} \in \Sigma \left[ s \right ]_2+
\Sigma\big[ \ve s \ve \big]_4 + \ve s \ve  \Sigma\big[ \ve s \ve \big]_2,
\]
which implies $\mu^{*} \le  \gamma^*$.

Second, consider the case $\| s^* \| > 0$,
then $\bt + 3 \sig \| s^* \| \ge 0$. Let $w = s - s^*$.
As in the proof of Theorem~2.2 of \cite{zhu2025CQR},
we have the expansion (note $B(s^*) s^* + g = 0$),
\begin{eqnarray}\label{eq:expansion}
M(s) -\mu^{*} & = & [B(s^*) s^* + g]^T w + \frac{1}{2} w^T B(s^*) w + \mc{F}_2  \\
& = &  \frac{1}{2} w^T B(s^*) w + \mc{F}_2,
\end{eqnarray}
where
\begin{eqnarray*}
\mc{F}_2 & = & \frac{1}{2} (\|s\| - \|s^*\| )^2
\Big[ \frac{\bt}{6} ( \|s^*\|+2\|s\|) + \frac{\sig}{2}(\|s^*\|+\|s\|)^2 \Big]  \\
&=&  \frac{1}{2} (\|s\| - \|s^*\|)^2 \Big[ \frac{\bt + 3 \sig \|s^*\|}{6}
    (\|s^*\| + 2 \|s\|)  +   \frac{\sig}{2}\|s\|^2 \Big].
\end{eqnarray*}
Since $\bt + 3 \sig \| s^* \| \ge 0$, we have
\[
\mc{F}_2 \in  \Sigma\big[ \ve s \ve \big]_4 + \ve s \ve  \Sigma\big[ \ve s \ve \big]_2.
\]
Since $s^*$ is a global minimizer, we have $B(s^*) \succeq 0$
by Theorem~\ref{thm:opcd}, so
\[
M(s) - \mu^{*} \in \Sigma \left[ s \right ]_2+
\Sigma\big[ \ve s \ve \big]_4 + \ve s \ve  \Sigma\big[ \ve s \ve \big]_2.
\]
This also implies $\mu^{*} \le  \gamma^*$.

In either case $\|s^*\| = 0$ or $\| s^* \| > 0$, we have shown
$\mu^{*} \le  \gamma^*$. On the other hand, it holds
$
\gamma^*  = \vartheta^* \le \mu^{*}
$
by Theorem~\ref{thm:bound}.
Therefore, we can further get $\gamma^* = \vartheta^* =  \mu^{*}$, i.e.,
the SDP relaxations  \reff{problem:dual:sos} and \reff{pro-moment} are tight.

\medskip \noindent
``$\Rightarrow$"  Assume the SDP relaxations \reff{problem:dual:sos} and \reff{pro-moment}
are tight. We need to show that \reff{bt>=-3} holds.
Note that if $s^* = 0$, then \reff{bt>=-3} automatically holds.
So we consider the general case that $s^* \ne 0$.
Since the relaxation \reff{problem:dual:sos} is tight,
the optimal values $\gamma^*$, $\mu^*$ and  $M(s^*)$ are equal
(i.e., $\gamma^* = \mu^* = M(s^*)$),
which is achievable for \reff{problem:dual:sos} by Theorem~\ref{thm:bound}.
So,
\[
M(s) - M(s^*) = [s]_1^T Q [s]_1 + \rho( \|s\| ),
\]
where $Q$ is a symmetric psd matrix and
$\rho( \| s \| ) \in \Sigma\big[ \ve s \ve \big]_4 +
\| s \| \Sigma\big[ \ve s \ve \big]_2$.
Note that $\rho$ is a univariate polynomial in $\| s \|$
and $\rho(z) \ge 0$ for all $z \ge 0$.
By Theorem~3.2.1 of \cite{mpo}, we get
\[
\rho( \|s \| ) \,  = \,  p( \|s\| )^2 + \|s\| q(\|s\|)^2
\]
for a quadratic polynomial $p(z)$ and a linear polynomial $q(z)$, hence
\[
M(s) - M(s^*) = [s]_1^T Q [s]_1 + p( \|s\| )^2 + \|s\| q(\|s\|)^2.
\]
Note that
\[
0 \, = \, [s^*]_1^T Q [s^*]_1 + p( \|s^*\| )^2 + \|s^*\| q(\|s^*\|)^2 .
\]
Because all the terms in the right hand side above are nonnegative,
\[
[s^*]_1^T Q [s^*]_1 = p( \|s^*\| )^2  =  \|s^*\| q(\|s^*\|)^2 = 0 .
\]
Since $s^* \ne 0$, we have $\| s^* \| > 0$ and hence
\[
 p( \|s^*\| )  =   q(\|s^*\|)  =  0.
\]
Without loss of generality, assume the leading coefficients of $p$ and $q$ are nonnegative.
The coefficient of $\| s \|^4$ of $p(\| s \|)^2$ is $\frac{\sig}{4}$,
so we can factorize $p,q$ as
\begin{eqnarray*}
p( \|s\| )  & = &  \frac{ \sqrt{\sig} }{2} (\| s \| - \| s^*\|) (\| s \| + \tau ), \\
q( \|s\| )  & = &  \sqrt{\eta} (\| s \| - \|s^*\|),
\end{eqnarray*}
for some real scalars $\tau \in \re$ and $\eta \ge 0$.
Hence, it holds
\[
 p( \|s\| )^2 + \|s\| q(\|s\|)^2 = (\| s \| - \|s^*\|)^2
\big[ \frac{\sig}{4} (\| s \| + \tau )^2 + \eta \|s\|  \big]
\]
\[
= \Big(\| s \|^2 - 2 \|s \| \|s^*\| + \| s^* \|^2 \Big)
\Big[ \frac{\sig}{4} \| s \|^2  + ( \frac{\sig \tau }{2} +
\eta) \|s\| + \frac{\sig\tau^2}{4} \Big] .
\]
To match the representation of $M(s) - M(s^*)$,
the coefficient of $\| s \|$ in $ p( \|s\| )^2 + \|s\| q(\|s\|)^2$
must be zero, so
\[
\|s^*\|^2 ( \frac{\sig \tau}{2} + \eta )  - \|s^*\|\frac{\sigma \tau^2}{2} =  0.
\]
The coefficient of $\| s \|^3$  in $ p( \|s\| )^2 + \|s\| q(\|s\|)^2$
must be $\frac{\beta}{6}$, so
\[
( \frac{\sig \tau}{2} + \eta )  - \frac{\| s^* \| \sig }{2}   =  \frac{\bt}{6}.
\]
The above equations imply that
\[
\|s^*\|^2 ( \frac{\sig \tau}{2} + \eta )  = \| s^* \| \frac{\sigma \tau^2}{2}, \quad
( \frac{\sig \tau}{2} + \eta )  = \frac{1}{6} ( \bt + 3 \sig \|s^*\| ).
\]
Therefore, it holds
\[
  \frac{\|s^*\|^2}{6} ( \bt + 3 \sig \|s^*\| )  \, =  \,
  \|s^*\|  \frac{\sigma \tau^2}{2} \ge 0.
\]
Since $\| s^* \| \ge  0$ and $\sigma \ge 0$,
the condition \reff{bt>=-3} must hold.

If \reff{bt>=-3} holds for one minimizer $s^*$, then
the relaxations \reff{problem:dual:sos} and \reff{pro-moment} are tight.
This further implies that \reff{bt>=-3} must also hold for
all other minimizers, if they exist.
Therefore, if \reff{bt>=-3} holds for one minimizer,
then it holds for all minimizers.
\qed
\end{proof}

\begin{rmk}\label{rmk:tight_cond}
The following can be implied by Theorem~\ref{thm:tight_cond}.
\begin{itemize}

\item [(i)]
If $s^*$ is a nonzero global minimizer
of the CQR problem (\ref{CQR:POLY}) and the SDP relaxations
\reff{problem:dual:sos} and \reff{pro-moment} are tight,
then $\bt + 3 \sig \| s^* \| \ge 0$,
which is the sufficient global optimality condition given in
\cite[Theorem~2.2]{zhu2025CQR}
(see also Theorem~\ref{thm:opcd2} in Section~\ref{sc:pre}).
This means that the relaxations
\reff{problem:dual:sos} and \reff{pro-moment} are tight if and only if
all nonzero global minimizers of (\ref{CQR:POLY})
satisfy this sufficient global optimality condition.

\item [(ii)]
For polynomial optimization problems, the classical Moment-SOS relaxation
is tight for the lowest relaxation order, under the {\it SOS-convexity}\footnote{A polynomial $p$ is said to be SOS-convex if there exists a polynomial matrix $H$ such that $\nabla^2 p = H^TH$.
SOS-convex polynomials must be convex.} assumption \cite[Sec.~7.2]{mpo}.
For the CQR problem (\ref{CQR:POLY}), we do not assume the objective $M(s)$ is convex or SOS-convex.
Indeed, if the parameter $\beta \ne 0$ (i.e., the third-order regularization term exists),
the objective $M(s)$ is not a polynomial, and it cannot be SOS-convex.
Moreover, when $H$ is not positive semidefinite, the objective $M(s)$ cannot be convex,
but the SDP relaxations \reff{problem:dual:sos} and \reff{pro-moment} are tight, as shown in Corollary~\ref{cor:eigH<0}.

\item [(iii)]
In general, tightness of the SDP relaxations (\ref{problem:dual:sos}) and (\ref{pro-moment})
is independent of the convexity of the CQR problem (\ref{CQR:POLY}).
That is, the relaxations may, or may not, be tight for convex CQR problems.
When $H$ is positive semidefinite and $\beta \ge 0$,
the CQR problem is clearly convex and the SDP relaxations (\ref{problem:dual:sos}) and (\ref{pro-moment}) are tight.
On the other hand, even if $M(s)$ is convex, the relaxation (\ref{pro-moment})
may not be tight. For instance, this is the case for
\be \label{oldexm:nottight:M(s)}
M(s) \, = \, s^4 -4 |s| ^3  + 6 s ^2-4 s  +1, \quad s \in \re^1.
\ee
One can verify that $M(s)$ is convex,
the minimum value $\mu^*=0$ but $\vartheta^* = -1$,
so the relaxation (\ref{pro-moment}) is not tight.
We refer to \cite{mpo,Slot2025} for relevant work on convex polynomial optimization.

\item [(iv)]
When $g^Tv = 0$ for the eigenvector $v$ of $H$ associated to its smallest eigenvalue $\lmd_{min}$, it may be difficult to find global minimizers for CQR problems
by solving its secular equation, particularly when $B(s^*)$ is singular.
This is usually referred to as the ``hard-case" for CQR and trust-region problems; see \cite{cartis2022evaluation,zhu2025CQR}.
For the SDP relaxations (\ref{problem:dual:sos}) and (\ref{pro-moment}), the condition (\ref{bt>=-3}) for tightness is in general independent of $g^Tv = 0$.
The relaxations \reff{problem:dual:sos} and \reff{pro-moment} may still be tight
for this hard case, even if $B(s^*)$ is singular.
We refer to Examples~\ref{ex:cases} (iv)-(v),
\ref{ex:two_isolated}, \ref{ex:6-1}, \ref{ex:sec6:add:ug=0},
\ref{ex:sec6:add:ug=0:pos} and \ref{ex:ar3:com:n=5} for such expositions.

\end{itemize}
\end{rmk}

When does the condition \reff{bt>=-3} hold? An interesting case is that $H$
has a nonpositive eigenvalue. Since $s^*$ is a global minimizer,
the expansion formula (\ref{eq:expansion}) shows
\be\label{eq:expansion>=0}
\frac{1}{2} w^T B(s^*) w + \mc{F}_2  \ge  0
 \quad \text{for all} \quad  w \in \re^n.
\ee
By Theorem~\ref{thm:opcd},
the global optimality condition implies
\be\label{eq:Bs}
B(s^*) = H + \frac{\bt}{2}  \|s^*\| I_n    + \sig   \|s^*\|^2 I_n \succeq 0.
\ee
Suppose $H$ has an eigenvalue $\lmd \le 0$.
Then, the above implies
\[
\frac{\bt}{2}   \|s^*\|  + \sig   \|s^*\|^2  \ge -\lmd \ge 0.
\]
Since $\sig \ge 0$, one can see that
\[
 \|s^*\| ( \bt  + 3 \sig   \|s^*\| )  \ge
 2 \big(\frac{\bt}{2}   \|s^*\|  + \sig   \|s^*\|^2 \big)  \ge  0.
\]
So, \reff{bt>=-3} holds, and we get the following corollary.

\begin{cor}  \label{cor:eigH<0}
If either $\bt \ge 0$ or $H$ has a nonpositive eigenvalue,
then the SDP relaxations \reff{problem:dual:sos} and \reff{pro-moment}
are tight, i.e., $\gamma^* = \vartheta^* =  \mu^*$.
\end{cor}

\begin{rmk}
The following can be implied by Corollary~\ref{cor:eigH<0}.

\begin{itemize}
\item [(i)] For the ARC problem (\ref{cubic-newton}),
we have $\sigma = 0$ and $\beta > 0$.
By Corollary~\ref{cor:eigH<0}, the SDP relaxations
\reff{problem:dual:sos} and \reff{pro-moment}
are tight for solving \reff{cubic-newton}.

\item [(ii)]
In \cite[Theorem~2.5]{zhu2025CQR}, there are five cases such that a root to a univariate  equation (with the specifically selected branch;
see \cite[Sections~2.2.2 and 2.2.3]{zhu2025CQR}) gives a global minimizer.
Case~1 of \cite[Theorem~2.5]{zhu2025CQR} is that $H$ is indefinite,
and Case~2 is that $H$ is positive definite and $\beta \ge 0$.
Interestingly, the relaxations \reff{problem:dual:sos} and \reff{pro-moment}
are tight for these two cases, by Corollary~\ref{cor:eigH<0}.
Moreover, let $\lmd_{min}$ be the smallest eigenvalue of $H$
and suppose $g^Tv \ne 0$ for all eigenvectors associated to $\lmd_{min}$.
For Case~4 that $H\succ 0$ and $\beta\le -3\sqrt{2\lmd_{min}\sigma}$,
there exists a global minimizer $s^*$ for (\ref{CQR:POLY})
such that $\beta + 3\sigma\Vert s^*\Vert \ge 0$.
This is shown in the proof for Case~4 of Theorem~2.5 in \cite{zhu2025CQR}.
Thus, by Theorem~\ref{thm:tight_cond}, the relaxations \reff{problem:dual:sos} and \reff{pro-moment} are tight for this case. For the other two cases,
they may or may not be tight (see Example~\ref{ex:cases}).
\end{itemize}

\end{rmk}

Theorem~\ref{thm:tight_cond} implies that
the SDP relaxations \reff{problem:dual:sos} and \reff{pro-moment}
are tight if and only if $\| s^* \| (\bt  + 3 \sig \| s^* \|) \ge 0$
for every minimizer $s^*$ of \reff{CQR:POLY}.
Indeed, if this inequality holds for one minimizer,
then it holds for all minimizers. Moreover, when the relaxations are tight,
we can show that all nonzero global minimizers have the same length.

\begin{thm}\label{thm:uniquenorm}
Assume at least one of $\sigma$ and $\beta$ is nonzero and
the SDP relaxations \reff{problem:dual:sos} and \reff{pro-moment} are tight.
If $s^*,\hat{s}$ are two minimizers of \reff{CQR:POLY}, then
\be  \label{samenorm}
 \| \hat{s} \|  \cdot  \| s^* \| \cdot \big( \| \hat{s} \| -  \| s^* \| \big) \, = \, 0.
\ee
\end{thm}
\begin{proof}
If one of $\| \hat{s} \|$ and $\| s^* \|$ is zero, then \reff{samenorm}
automatically holds. Now suppose they are both nonzero.
By Theorem~\ref{thm:bound}, the optimal value of \reff{problem:dual:sos} is achievable.
Since \reff{problem:dual:sos} and \reff{pro-moment} are tight relaxations,
as in the proof of Theorem~\ref{thm:tight_cond}, we can show that
there exist a psd matrix $Q$,
a quadratic polynomial $p(z)$ and a linear polynomial $q(z)$ such that
\[
M(s) - \mu^* = [s]_1^T Q [s]_1 + p( \|s\| )^2 + \|s\| q(\|s\|)^2.
\]
Since $\hat{s}$, $s^*$ are minimizers of (\ref{CQR:POLY}), one can see
\[
 M(s^*) - \mu^* = M(\hat{s}) - \mu^* = 0.
\]
Both $\| \hat{s} \|$ and $\| s^* \|$ are positive, so
\[ p( \|\hat{s}\| ) = p( \|s^*\|) = q(\|\hat{s}\|) = q(\|s^*\|) = 0. \]
We discuss in two cases.

\bit

\item Suppose $q$ is not identically zero. Then, both $\|\hat{s}\|$ and $\|s^*\|$
are roots of $q(z) = 0$. Since the degree of $q$ is one,
we can get $\|\hat{s}\|=\|s^*\|$.

\item Suppose $q$ is identically zero.
As in the proof of Theorem~\ref{thm:tight_cond}, we can factorize that
\[
p( \|s\| )   =  \frac{ \sqrt{\sig} }{2} (\| s \| - \| s^*\|)  (\| s \| + \tau )
\]
for some real scalar $\tau$.
Since $q = 0$ and the coefficient of $\Vert s\Vert^3$ (resp., $\Vert s\Vert$)
in $p( \|s\| )^2$ is $\frac{\beta}{6}$ (resp., $0$), we have
\[
\frac{\sig}{2} (\tau  - \| s^* \| )  =  \frac{\bt}{6}, \quad
\frac{\sig}{2} \tau \| s^* \| \big(
\|s^*\|    - \tau  \big) = 0 .
\]
Note $\|s^*\|>0$ and the above implies $\beta = 0$ if $\sigma = 0$.
Since $\sigma$ and $\beta$ are not both zero, we must have $\sig >0$.
If $\tau \ne 0$, then $\tau =\|s^*\|$ and hence
\[  p( \|s\| )   =  \frac{ \sqrt{\sig} }{2} (\| s \|^2 - \| s^*\|^2). \]
Since $p( \|\hat{s}\| ) = 0$, we get $\| \hat{s} \| = \| s^*\|$.
If $\tau = 0$, then
\[  p( \|s\| )   =  \frac{ \sqrt{\sig} }{2} \|s\| (\| s \| - \| s^*\|). \]
Since $p( \|\hat{s}\| ) = 0$ and $\| \hat{s} \| > 0$,
we also get $\| \hat{s} \| = \| s^*\|$.

\eit
Therefore, we get $\| \hat{s} \| = \| s^* \|$ for both cases.
\qed
\end{proof}

The following are some exposition examples for the above conclusions for the CQR problems.
\begin{Example}\label{ex:cases} \rm
\noindent
(i) Consider the CQR problem
\[
\min_{s\in \re^2} \quad  \Vert s\Vert ^4- \frac{8}{3}\Vert s \Vert ^3+2s_1^2+s_2^2-8s_1+1.
\]
In terms of the formulation of \eqref{CQR:POLY}, this corresponds to
\[ \sigma=4,\quad \beta=-16,\quad H=\bbm
4 &   0  \\
0  &   2 \\
\ebm,\quad
g=-\bbm
8  \\
0 \\
\ebm, \quad f_0=1.\]
By solving \reff{pro-moment}, we can get the optimizer
\[
Y^*=
\bbm
1  &  2 & 0  \\
2  &  4 & 0  \\
0  &  0 & 0  \\
\ebm, \quad
Z_1^*=
\bbm
1   &  2   &  4 \\
2   &  4   &  8  \\
4   &  8   &  16  \\
\ebm, \quad
Z_2^*=
\bbm
2 &   4  \\
4 &   8 \\
\ebm.
\]
The optimal values of \reff{problem:dual:sos} and \reff{pro-moment}
are $\gamma^* = \vartheta^* = -\frac{37}{3}$.
Moreover, one may check that $M(s^*) = -\frac{37}{3}$ when $s^* = (2,0)$,
which implies that $s^*$ is a global minimizer.
The SDP relaxations \reff{problem:dual:sos} and \reff{pro-moment} are tight, and condition~\reff{bt>=-3} holds. Besides that, since
\[
\beta < 0,\quad H\succ 0,\quad \beta < -3\sqrt{2\sigma\lmd_{min}} = -12,
\]
one may check that this problem belongs to Case~4 of \cite[Theorem~2.5]{zhu2025CQR}. \\
\noindent
(ii) Consider the CQR problem
\[
\min_{s\in \re^2} \quad  \Vert s\Vert ^4-6 \Vert s \Vert ^3+13s_1^2+13s_2^2+12s_1+4.
\]
This corresponds to parameters
\[ \sigma=4,\quad \beta=-36,\quad H=26I_2,\quad g=\left[\begin{array}{c}
 12 \\ 0
\end{array}\right],\quad f_0=4.\]
By solving \reff{pro-moment}, we get the optimizer
\[
Y^*=
\bbm
1  &   -1.5 & 0 \\
-1.5  &   2.25 & 0  \\
0 & 0 & 0
\ebm, \quad
Z_1^*=
\bbm
1   &   0.75  &   2.25  \\
0.75  &   2.25  &   6.75 \\
2.25 &   6.75  &  20.25  \\
\ebm, \quad
Z_2^*=
\bbm
0.75  &   2.25  \\
2.25  &   6.75 \\
\ebm.
\]
The optimal values of \reff{problem:dual:sos} and \reff{pro-moment}
are $\gamma^* = \vartheta^* =  -5$.
There are two global optimizers: $s^* = (-1,0)$ or $s^* = (-2,0)$.
The minimum value $\mu^* = 0$, so the SDP relaxations
\reff{problem:dual:sos} and \reff{pro-moment} are not tight.
Note $B(s^*) = 12I_2$ at $s^* = (-1,0)$, $B(s^*) = 6I_2$ at $s^* = (-2,0)$,
and $\| s^* \| >0$. The sufficient global optimality condition
$\bt + 3 \sig \|s^*\| \ge 0$ fails to hold.
These two global minimizers have different Euclidean norms.
Theorem~\ref{thm:uniquenorm} also implies
\reff{problem:dual:sos} and \reff{pro-moment} are not tight.
This problem belongs to Case~3 of \cite[Theorem~2.5]{zhu2025CQR}
that $\beta < 0$, $H\succ 0$, and
\[
-4\sqrt{\sigma \lmd_{min}}<\beta <0,\quad
-\frac{\beta}{4\sigma} >
\Big\Vert (H-\frac{\beta^2}{16\sigma} I_n )^{-1}g \Big\Vert.
\]
\noindent
(iii) Consider the CQR problem
\[
\min_{s\in \re^2} \quad  \Vert s\Vert ^4-4 \Vert s \Vert ^3+10s_1^2+5s_2^2-12s_1-2s_2+1.
\]
In terms of the formulation of \eqref{CQR:POLY}, this corresponds to
\[ \sigma=4,\quad \beta=-24,\quad H=\bbm
20 &   0  \\
0  &   10 \\
\ebm,\quad
g=\bbm
-12  \\
-2 \\
\ebm, \quad f_0=1.\]
By solving \reff{pro-moment}, we can get the optimizer
\be\label{eq:nontight_moment_minimizer}
Y^*=
\bbm
1  &  1 & 1  \\
1  &  1 & 1  \\
1  &  1 & 1  \\
\ebm, \quad
Z_1^*=
\bbm
1   &  1   &  2 \\
1   &  2   &  4  \\
2   &  4   &  8  \\
\ebm, \quad
Z_2^*=
\bbm
1 &   2  \\
2 &   4 \\
\ebm.
\ee
The optimal values of \reff{problem:dual:sos} and \reff{pro-moment}
are $\gamma^* = \vartheta^* = -6$.
However, the minimum value $\mu^* \approx -5.6421$,
which can be obtained by solving the equivalent polynomial optimization problem
\reff{eq:trans_pop} using the Moment-SOS relaxation methods \cite{mpo}.
The SDP relaxations \reff{problem:dual:sos} and \reff{pro-moment} are not tight.
The condition~\reff{bt>=-3} fails to hold.
Besides that, one may check that this problem belongs to Case~5 of
\cite[Theorem~2.5]{zhu2025CQR} that $\beta < 0$, $H\succ 0$, and
\[
-4\sqrt{\sigma \lmd_{min}}<\beta <0,\quad
-\frac{\beta}{4\sigma} < \Big\Vert
(H-\frac{\beta^2}{16\sigma} I_n )^{-1}g \Big\Vert.
\]
\noindent
(iv) When $s^* = 0$ is a global minimizer of CQR problem,
the condition $\bt \ge -3 \sig \| s^* \|$ is not necessary for
\reff{problem:dual:sos} and \reff{pro-moment} to be tight.
Consider the CQR problem
\[
\min_{s\in \re^n} \quad  \Vert s\Vert^2 (\| s \| - 1 )^2   + \eps \| s \|^3 ,
\]
for a parameter $\eps \in (0, 2)$. The corresponding parameters are
\[
\sigma=4,\quad \beta=-12 + 6\eps,
\quad H = 2I_n,\quad g =  \mathbf{0}_n,  \quad f_0 = 0.
\]
The minimum value $\mu^* = 0$ and the unique minimizer is $s^* = \mathbf{0}_n$.
The relaxations \reff{problem:dual:sos} and \reff{pro-moment} are tight, i.e.,
$\gamma^* = \vartheta^* =  0$, since
\[
\Vert s\Vert^2 (\| s \| - 1 )^2   + \eps \| s \|^3 - 0
\in \Sigma[\Vert s\Vert]_4 + \Vert s\Vert\Sigma[\Vert s\Vert]_2.
\]
Clearly, condition~\reff{bt>=-3} holds.
We remark that since $g = \mathbf{0}_n$, one has $g^Tv = 0$
holds for all eigenvector $v$ of $H$ corresponds to $\lmd_{min}$.
Moreover, one may check that when
$\epsilon > 2-\frac{4\sqrt{2}}{3} \approx  0.1144$,
this problem belongs to Case~3 of \cite[Theorem~2.5]{zhu2025CQR};
when $0 < \epsilon < 2-\frac{4\sqrt{2}}{3}$, it belongs to Case~5.
For $\epsilon = 2-\frac{4\sqrt{2}}{3}$,
the matrix $H-\frac{\beta^2}{16\sigma} I_n$ is not invertible.\\
\noindent
(v) Consider the CQR problem in Example~\ref{ex:rank1_inf_global_min}.
The corresponding parameters are:
\[
\sigma=4,\quad \beta=-24,\quad  H = 8I_n, \quad g =  \mathbf{0}_n,
\quad f_0 = 0.
\]
The minimum value $\mu^* = 0$, and the global minimizers are
$0$ and all points $s^*$ such that $\| s^* \| = 2$.
As shown in Example~\ref{ex:rank1_inf_global_min}, the relaxations
\reff{problem:dual:sos} and \reff{pro-moment} are tight.
One can check that the condition~\reff{bt>=-3} holds.
This problem belongs to Case~4 of \cite[Theorem~2.5]{zhu2025CQR}
since $H\succ 0$ and $\beta = -3\sqrt{2\sigma \lmd_{min}}$.
For all nonzero minimizers $s^*$, it holds
\[
B(s^*) = 8I_n+\frac{\beta}{2}\Vert s^*\Vert I_n + \sigma \Vert s^*\Vert^2 I_n
= {\bf 0}_{n\times n}.
\]
The matrix $B(s^*)$ is not invertible.
\end{Example}

\section{Certifying tightness and extracting minimizers}
\label{sc:alg}

This section discusses how to detect tightness of the relaxations
\reff{problem:dual:sos} and \reff{pro-moment},
and how to extract global minimizers of \eqref{CQR:POLY}.

As shown in Theorem~\ref{thm:bound},
the optimal value $\gamma^*$ of \reff{problem:dual:sos} is attainable.
There exist symmetric matrices $X_0^*, X_1^*, X_2^*$ such that
\[
M(s) -\gamma^* \, = \,  \bbm 1 \\ s \ebm^T X_0^* \bbm 1 \\ s \ebm +
\bbm 1 \\ \| s \| \\ \| s \|^2 \ebm^T X_1^* \bbm 1 \\ \| s \| \\ \| s \|^2 \ebm  +
\|s \| \bbm 1 \\ \| s \|  \ebm^T X_2^* \bbm 1 \\ \| s \|   \ebm ,
\]
\[
X_0^* \in \mc{S}_+^{n+1}, \quad
X_1^* \in \mc{S}_+^{3}, \quad
X_2^* \in \mc{S}_+^{2}.
\]
Suppose $(Y^*, Z_1^*, Z_2^*)$ is an optimizer of \reff{pro-moment}.
Let $\ell^*$ be the linear functional determined by $(Y^*, Z_1^*, Z_2^*)$
as in \reff{mat:Y}-\reff{mat:Z2}.
By Theorem~\ref{thm:bound}, there is no duality gap between
\reff{problem:dual:sos} and \reff{pro-moment}, thus
\[
0 = \vartheta^* - \gamma^* = \vartheta(Y^*, Z_1^*, Z_2^*) - \gamma^*
= \ell^* (M(s) -\gamma^*) =
\]
\[
Y^* \bullet X_0^* + Z_1^* \bullet X_1^* +  Z_2^* \bullet X_2^* .
\]
Since all the matrices are psd, the above implies that
\[
Y^*X_0^* = {\bf 0}_{(n+1)\times (n+1)},\quad
Z_1^* X_1^* = {\bf 0}_{3\times 3},\quad
Z_2^* X_2^* = {\bf 0}_{2\times 2}.
\]
The first equation in the above implies that the column space of $X_0^*$
belongs to the null space of $Y^*$, whose dimension is $n+1-\rank \, Y^*$.
Since $\rank \, X_0^*$ equals the dimension of the column space of $X_0^*$, it holds
\[
\rank \, Y^* + \rank \, X_0^* \le n+1.
\]
Similarly, we can get
\[
\rank \,  Z_1^* + \rank \,  X_1^* \le 3, \quad
\rank \,  Z_2^* + \rank \,  X_2^* \le 2 .
\]
We can write $X_0^*=R^T R$ for some $R\in \re^{ r\times (n+1)}$ with $r\le n+1$.
Note $\rank \,  Z_1^* \ge 1$, so $r_1 \coloneqq \rank \,  X_1^* \le 2$.
Also note $r_2 \coloneqq \rank \,  X_2^* \le 2$.
Then, there exist $a_1,a_2 \in \re^{3}$ and
$b_1, b_2 \in \re^{2}$ such that
\[
  X_1^* \, = \, a_1a_1^T+ a_2 a_2^T, \quad
  X_2^* \, = \, b_1b_1^T + b_2b_2^T.
\]
In the above, we let $a_2 = 0$ if $r_1 = 1$, $a_1 = a_2 = 0$ if $r_1 = 0$,
$b_2 = 0$ if $r_2 = 1$, and $b_1 = b_2 = 0$ if $r_2 = 0$.
We assume that $\sigma$ and $\beta$ are not simultaneously zero,
so at least one of $r_1$ and $r_2$ must be nonzero. Then
\[
\begin{array}{ccl}
\bbm 1 \\ \| s \| \\ \| s \|^2 \ebm^T X_1^* \bbm 1 \\ \| s \| \\ \| s \|^2 \ebm
& = &\bbm 1 \\ \| s \| \\ \| s \|^2 \ebm^T
\big( a_1a_1^T  +  a_2a_2^T  \big) \bbm 1 \\ \| s \| \\ \| s \|^2 \ebm  \\
& = &  (a_{10}+a_{11}\ve s\ve +a_{12}\ve s\ve^2)^2 \\
&& \quad +(a_{20}+a_{21}\ve s\ve +a_{22}\ve s\ve^2)^2 ,
\end{array}
\]
\[
\begin{array}{ccl}
\bbm 1 \\ \| s \|  \ebm^T X_2^* \bbm 1 \\ \| s \|  \ebm
& = &\bbm 1 \\ \| s \|  \ebm^T
\Big(  b_1b_1^T  +  b_2b_2^T \Big) \bbm 1 \\ \| s \|  \ebm \\
& = & (b_{10}+b_{11}\ve s\ve )^2 + (b_{20}+b_{21}\ve s\ve )^2.
\end{array}
\]
In the above, $a_i = (a_{i0}, a_{i1}, a_{i2})$
and $b_i = (b_{i0}, b_{i1})$ for $i=1,2$.
Denote the univariate polynomials
\begin{eqnarray}
p_i(z) &=& a_{i0}+a_{i1}z + a_{i2}z^2,    \quad i=1,2,     \\
q_j(z) &=& b_{j0}+b_{j1}z,  \quad j=1, 2 .
\end{eqnarray}
Then, it holds
\be  \label{eq:Mc-gm}
M(s) - \gamma^* = [s]_1^TR^TR[s]_1 + \sum_{i=1}^{2}p_i(\Vert s \Vert )^2
+ \Vert s \Vert \sum _{j=1}^{2}q_j(\Vert s \Vert )^2.
\ee
Suppose $s^*$ is a global minimizer of \eqref{CQR:POLY}, then
\[
M(s^*) - \gamma^* \, = \, [s^*]_1^TR^TR[s^*]_1 +
\sum_{i=1}^{2}p_i(\Vert s^* \Vert )^2 +
\Vert s^* \Vert \sum _{j=1}^{2}q_j(\Vert s^* \Vert )^2.
\]
If the relaxation~\reff{problem:dual:sos} is tight,
i.e., $\gamma^* = M(s^*)$, then $s^*$ is a solution to
\be  \label{eq:3terms=0}
\left \{
\begin{array}{ll}
 R [s]_1  = 0,  \\
 p_1(\Vert s \Vert ) = p_2(\Vert s \Vert )  = 0, \\
 \Vert s \Vert q_1(\Vert s \Vert ) = \Vert s \Vert q_{2}(\Vert s \Vert ) = 0.
\end{array}
\right.
\ee
Conversely, if $s^*\in \re^n$ satisfies (\ref{eq:3terms=0}), then
\[
M(s^*) - \gamma^* = 0.
\]
Since $\gamma^* = \mu^* \le M(s^*)$,
the above implies $s^*$ is a global minimizer of \eqref{CQR:POLY}
and the relaxation \reff{problem:dual:sos} is tight.
The latter further implies \reff{pro-moment} is also tight.
Moreover, by Theorems~\ref{thm:tight_cond} and \ref{thm:uniquenorm},
there exists $z^*>0$ with $\bt  \ge  -3 \sig z^*$ such that
for all nonzero solutions $s^*$ to (\ref{eq:3terms=0}), if there are any, it holds
\[  \Vert s^*\Vert \, = \, z^*. \]

Summarizing the above, we get the following algorithm for
checking tightness of \reff{problem:dual:sos}
and extracting global minimizers for (\ref{CQR:POLY}).

\begin{alg} \label{alg:finds}
Let $\mathcal{S} \coloneqq  \emptyset$ and do the following:
\begin{itemize}[leftmargin = 5em]

\item [Step~1.]
Solve the SDP relaxations \reff{problem:dual:sos} and \reff{pro-moment},
get the optimal values $\gamma^*$, $\vartheta^*$,
and the representation \reff{eq:Mc-gm}.

\item [Step~2.]
If the zero vector $0$ satisfies \reff{eq:3terms=0},
let $\mathcal{S} \coloneqq \mathcal{S}\cup \{ 0 \}$.

\item [Step~3.] Solve the following system of univariate equations:
\[
p_1(z) = p_{2}(z) = q_1(z) = q_2(z) = 0.
\]
If they have a positive real zero $z^*$, then go to Step~4;
otherwise, output the set $\mc{S}$.

\item [Step~4.]  Output the following set and stop:
\[
\mc{S} \coloneqq \mc{S} \cup  \{ s\in \re^n: \ve s\ve =z^*,\, R [s]_1 = 0 \} .
\]

\end{itemize}
\end{alg}

In the above, we have seen that the system \reff{eq:3terms=0} has a solution
if and only if the relaxations \reff{problem:dual:sos} and \reff{pro-moment} are tight.
For such cases, we can get all minimizers of \reff{CQR:POLY}
by Algorithm~\ref{alg:finds}.
Thus, we get the following theorem.

\begin{thm}\label{th:tightness}
The SDP relaxations \reff{problem:dual:sos} and \reff{pro-moment}
are tight if and only if
the system (\ref{eq:3terms=0}) has a solution.
If they are tight, the set $\mc{S}$ output by Algorithm~\ref{alg:finds}
consists of all global minimizers of \reff{CQR:POLY};
if otherwise they are not tight, the output set $\mc{S}$ is empty.
\end{thm}

By Algorithm~\ref{alg:finds} and Theorem~\ref{th:tightness},
when the SDP relaxations \reff{problem:dual:sos} and \reff{pro-moment} are tight,
the set of minimizers for the CQR problem \reff{CQR:POLY}
is equal to the solution set of (\ref{eq:3terms=0}).
Note that at least one of $p_1$ and $q_1$ is nonzero. Denote
\[
L  \, \coloneqq \,  \{ s\in\re^n: R[s]_1 = 0  \}.
\]
Assume \reff{problem:dual:sos} and \reff{pro-moment}
are tight and \reff{CQR:POLY} has nonzero minimizers.
Then, there exists a unique $z^*$ such that
$p_1(z^*) = p_{2}(z^*) = q_1(z^*) = q_2(z^*) = 0,$
and the set of all nonzero minimizers is
\[ L \cap \{s\in\re^n: \Vert s \Vert = z^*\}. \]
Therefore, if \reff{CQR:POLY} has infinitely many minimizers,
then the dimension of $L$ must be greater than $1$.
For instance, consider the CQR problem (\ref{eq:s=0ors=2}) in Example~\ref{ex:rank1_inf_global_min}.
The optimal value for the SDP relaxation (\ref{problem:dual:sos}) is $\gm^* = 0$.
By (\ref{eq:rank1_inf_global_min_sos}),
the maximizer of (\ref{problem:dual:sos})
satisfies \reff{eq:3terms=0} with
\[
R = 0,\quad
p_1(\Vert s\Vert) = \Vert s\Vert^2-2\Vert s\Vert,\quad
p_2 = q_1=q_2 = 0.
\]
Algorithm~\ref{alg:finds} produces the set of all
global minimizers for (\ref{eq:s=0ors=2}), which is
\[\mc{S} = \{0\} \cup \{ s\in\re^n: \Vert s\Vert = 2\} .\]
The following two examples are expositions for Algorithm~\ref{alg:finds}.

\begin{Example}\label{ex:two_isolated}  \rm
Consider the CQR problem
\[
\min_{s\in \re^{10}} \quad  7s_1^2-\sum_{i=2}^{10}s_i^2
+\sum_{i=1}^{5} 2 s_1s_{2i}-\sum_{i=1}^{4} 2  s_1s_{2i+1}+\ve s \ve ^4.
\]
The corresponding parameters are $\sigma=4$, $\beta=0$, and $g = \mathbf{0}_{10}$.
Moreover, the matrix $H$ is not positive definite.
By solving \reff{pro-moment}, we get the optimal triple $(Y^*, Z_1^*, Z_2^*)$
with $\rank\, Y^* = 2$, $\rank\, Z_1^* =\rank\, Z_2^* = 1$.
It took around $0.03$ second. Algorithm~\ref{alg:finds} produces
\[
\begin{array}{ccl}
  p_{1}(z) & \approx &-0.7071+0.7071z^2, \\
  p_{2}(z) & \approx  & 0.4082-0.8165z+0.4083z^2,\\
  q_{1}(z) & \approx  &-0.7071+0.7071z,
\end{array}
\]
and $R\in\re^{9\times 11}$, which has rank $9$.
For neatness, we do not display $R$ here.
Solving $p_1(z)=p_2(z)=q_1(z)=0$,
we can get $z^* = 1$.
In Step~4, the minimizers are given as
\[
\begin{array}{ccl}
 s^* &=&  t\, (-1,1,-1,1,-1,1,-1,1,-1,1), \\
 \Vert s^* \Vert   &=& 1, \quad t\in \re.
\end{array}
\]
Therefore, we get two global minimizers
\[
s^* \,  = \, \frac{\pm 1}{\sqrt{10}}  (-1,1,-1,1,-1,1,-1,1,-1,1).
\]
The optimal values are $\gamma^* = \vartheta^* =  \mu^* = -1$.
The relaxations \reff{problem:dual:sos} and \reff{pro-moment} are tight.
\end{Example}

Theorem~\ref{th:tightness} shows that the SDP relaxations
\reff{problem:dual:sos} and \reff{pro-moment} are tight
if and only if Algorithm~\ref{alg:finds} produces a nonempty set of solutions $\mc{S}$.
In other words, when \reff{problem:dual:sos} and \reff{pro-moment}
are not tight, the nontightness can also be detected by Algorithm~\ref{alg:finds}.
The following is such an exposition.

\begin{Example} \rm \label{sec5:exm5.4}
Consider the CQR problem in Example~\ref{ex:cases}(iii).
The optimal values of the SDP relaxations \reff{problem:dual:sos} and \reff{pro-moment}
are $\gm^*=\vartheta^* = -6$, and the minimizer matrices of (\ref{pro-moment})
are given by (\ref{eq:nontight_moment_minimizer}).
Also, for the optimal solutions of (\ref{problem:dual:sos}), we have
$\rank\, X_0^* = 2$, $\rank\, X_1^* = 1$, $\rank\, X_2^* = 1$.
Algorithm~\ref{alg:finds} produces
\[
\begin{gathered}
R = \left[\begin{array}{rrr}
  -0.7368 &  0.6731 &  0.0637\\
  -0.3518 & -0.4622 &  0.8140
  \end{array}
\right], \\
p_1(z) = -0.0001 + 0.8944z -0.4472z^2, \\
q_1(z) = -0.8944 + 0.4473z.
\end{gathered}
\]
One can check that (\ref{eq:3terms=0}) does not have a solution.
This shows that the relaxations \reff{problem:dual:sos} and \reff{pro-moment}
are not tight, by Theorem~\ref{th:tightness}.
\end{Example}

\begin{rmk}
When the SDP relaxations \reff{problem:dual:sos} and \reff{pro-moment}
are not tight, one can get an approximate solution from the matrix $Y^*$
by setting $s^* = (Y^*_{10}, Y^*_{20}, \ldots, Y^*_{n0})$.
For Example~\ref{sec5:exm5.4}, solving the SDP relaxation \reff{pro-moment}
gives $\vartheta^* = -6$ and the optimizer matrix
\[
Y^* = \begin{bmatrix}
1 & 1 & 1 \\
1 & 1 & 1 \\
1 & 1 & 1
\end{bmatrix}.
\]
For $s^* = (1, 1)$, the objective value
is $M(s^*) \approx -5.3137$,
which is an upper bound for $\mu^*$.
When \reff{problem:dual:sos} and \reff{pro-moment}
are not tight, how to choose the best approximate minimizer
is mostly open for the authors.
The following guidelines provide practical strategies for solving CQR problems when the relaxations \reff{problem:dual:sos} and \reff{pro-moment} are detected to be not tight.
\begin{itemize}

\item One may refer to the secular equation approach
(e.g., the Newton-Cholesky method in \cite{zhu2025CQR}).
Indeed, it can be regarded as the default method in regular cases when the CQR problem is large scale, while the SDP relaxations \reff{problem:dual:sos} and \reff{pro-moment} can provide certificates of global optimality and are particularly useful in hard cases (see Remark~4.2 (iv)).

\item For CQR problems appearing in the scheme of adaptive regularization methods,
if the SDP relaxations are not tight, then one may consider updating the parameters,
for example, increasing the value of $\sigma$ of the CQR problem.
Although this will change the CQR problem, it is fully consistent
with regularization based outer algorithms.
Indeed, increasing $\sigma$ is a natural strategy for adaptive regularization methods.
We refer to \cite{Zhu2026Sufficiently} for related work on AR$3$ problems.
\end{itemize}
\end{rmk}

\begin{rmk}
We make some discussions about the number of
nonzero global minimizers of (\ref{CQR:POLY}).
When \reff{problem:dual:sos} and \reff{pro-moment} are tight,
let $z^* >0$ be the length as in Step~3 of Algorithm~\ref{alg:finds}.
By Theorems~\ref{thm:uniquenorm} and \ref{th:tightness},
a point $s\in\re^n$ is a nonzero global minimizer of (\ref{CQR:POLY}) if and only if
\be  \label{eq:nonzero_sol_iff}
\ve s\ve^2 =(z^*)^2,\quad R [s]_1 = 0.
\ee
Let $\hat{\mc{S}}$ be the solution set to \reff{eq:nonzero_sol_iff}.
Then $\hat{\mc{S}}$ is the set of nonzero global minimizers of (\ref{CQR:POLY}).
Note that the degree of \reff{eq:nonzero_sol_iff} is two.
By quadratic formula, if \reff{eq:nonzero_sol_iff} has only finitely many solutions,
then the cardinality of $\hat{\mc{S}}$ is at most two.
Therefore, if (\ref{CQR:POLY}) has finitely many global minimizers
and the SDP relaxations \reff{problem:dual:sos} and \reff{pro-moment} are tight,
then there are at most two nonzero minimizers for (\ref{CQR:POLY}),
hence \reff{CQR:POLY} has at most three global minimizers
if the zero point is also one of them.
For instance, there are $3$ global minimizers for
$M(s) = s^2(|s|-1)^2$ in $s \in \re^1$.
\end{rmk}

\section{Numerical experiments}
\label{sc:num}

We present some numerical experiments for solving
the cubic-quartic regularization problem \reff{CQR:POLY}.
The SDP relaxations \reff{problem:dual:sos} and \reff{pro-moment}
are solved by the software {\tt Mosek} \cite{mosek}.
The computations are implemented in MATLAB R2022b
on a Lenovo Laptop with CPU@2.10GHz and RAM 16.0G.
For neatness of presentation, all computational results
are displayed in four decimal digits.

Recall that $\mu^*$ denotes the minimum value of \reff{CQR:POLY},
$\gamma^*$ denotes the maximum value of \reff{problem:dual:sos},
$\vartheta^*$ denotes the minimum value of \reff{pro-moment},
and $s^*$ denotes the computed global minimizer of \reff{CQR:POLY}.
We remark that by Theorem~\ref{thm:bound}, $\gamma^*$ gives a lower bound for $\mu^*$,
and $\mu^*\le M(s)$ for all $s\in\re^n$.
The SDP relaxations \reff{problem:dual:sos} and \reff{pro-moment}
are tight if $M(s^*) = \gamma^*$.
In computational practice, one typically cannot have
$M(s^*)$ equal to $\gamma^*$ exactly, due to numerical errors.
To measure the numerical accuracy of $s^*$,  we use the absolute and relative errors
\[
\texttt{err-abs}  =  \left |M(s^*)-\gamma^*\right |, \quad
\texttt{err-rel}  =  \left |\frac{M(s^*)-\gamma^*}{ M(s^*) }\right |.
\]

To check the efficiency of SDP relaxations \reff{problem:dual:sos} and \reff{pro-moment},
we compare with some existing methods for solving CQR problems.
There exists the Cholesky-Newton (Cho.-New.) method \cite[Algorithm 1]{zhu2025CQR}
by solving secular equations. For this method, we set the error tolerance to be $10^{-8}$
and apply the initialization trick described in Remark~2.1 of \cite{zhu2025CQR}.
We also apply the MATLAB function {\tt fminunc} with default parameter settings to solve the CQR problem directly.
Since {\tt Gurobi} \cite{Gurobi} does not support polynomials of degrees higher than two,
it cannot solve the CQR problem (\ref{CQR:POLY}) directly.
By introducing auxiliary variables $t_1$, $t_2$, $t_3$, the CQR problem (\ref{CQR:POLY})
can be equivalently reformulated as the quadratically constrained quadratic program (QCQP):
\be  \label{eq:trans_qcqp}
\left\{ \begin{array}{cl}
\min &
f_0+g^Ts+\frac{1}{2}s^THs+\frac{\beta}{6}t_2+\frac{\sigma}{4} t_3, \\
\st & t_1 = s^Ts,\ t_3 = t_1^2,\ t_2^2 = t_1t_3,\ t_2\ge 0 \\
       & s\in\re^n,\ t_1\in\re,\ t_2\in\re,\ t_3\in\re.
\end{array}
\right.
\ee
The software {\tt Gurobi} is then applied to solve \reff{eq:trans_qcqp}.
For CQR problems with the dimension $n \leq 10$,
we also apply the global optimization software {\tt Baron} \cite{Sahinidis1996}
to solve (\ref{CQR:POLY}) directly.
The current publicly available version of {\tt Baron}
does not solve problems of dimension bigger than ten.
To use {\tt Gurobi} and {\tt Baron}, we apply the MATLAB toolbox
{\tt Yalmip} \cite{yalmip} with the default parameter settings
and impose the additional box constraints $s\in [-10,10]^n$.
In the tables of computational results, ``SDP relax." stands for the relaxations (\ref{problem:dual:sos}) and (\ref{pro-moment}),
``$M_{value}$" represents the value of the function $M(s)$ at the point that is returned, and ``Time (s)" stands for the consumed time in seconds.

\begin{Example}\label{ex:6-1}  \rm
Consider the CQR problem
\[
\min_{s\in \re^5} \quad \sum_{1\le i<j\le 5}s_is_{j}-\frac{5}{2}\sum_{i=1}^{5}s_i^2
 - \ve s \ve ^3 + \ve s\ve^4.
\]
Solving \reff{problem:dual:sos} and \reff{pro-moment}, we get
$\rank\, Y^* = 5$, $\rank\, Z_1^* = \rank\, Z_2^* = 1$.
It took around $0.04$ second.
By Algorithm~\ref{alg:finds}, we get
\[
\begin{array}{ccl}
  p_{1}(z) & \approx &-0.9501+0.0703z+0.3040z^2,\\
  p_{2}(z) & \approx & 0.0926-0.8669z+0.4898z^2,\\
  q_{1}(z) & \approx &-0.8560+0.5170z, \\
  R & \approx & [\, 0,\, -0.4472 ,\,  -0.4472 ,\,  -0.4472 ,\,  -0.4472  ,\, -0.4472\, ].
\end{array}
\]
In Step~3, we get $z^* \approx  1.6559.$
The set of global minimizers can be parameterized as
$s^* = t_1\xi_1+t_2\xi_2+t_3\xi_3+t_4\xi_4$
with $\Vert s^*\Vert \approx 1.6559$, where
\[
\begin{array}{cclccl}
  \xi_1 &= & (1,\,-1,\,0,\,0,\, 0), \qquad
  \xi_2 &= & (1,\,0,\,-1,\,0,\, 0), \\
  \xi_3 &= & (1,\,0,\,0,\,-1,\, 0), \qquad
  \xi_4 &= & (1,\,0,\,0,\,0,\, -1). \\
\end{array}
\]
The optimal values are $\gamma^* = \vartheta^* =  \mu^* \approx -5.2479$.
The relaxations \reff{problem:dual:sos} and \reff{pro-moment} are tight.
In particular, $\hat{s}^* \approx (1.1709,\,-1.1709,\,0,\,0,\, 0)$
is a global minimizer, at which the errors are
\[
\begin{array}{c}
\texttt{err-abs}  =  1.08 \cdot 10^{-8}, \quad
\texttt{err-rel}  =  2.06 \cdot 10^{-9}.
\end{array}
\]
Performance of the Cholesky-Newton method, {\tt fminunc}, {\tt Baron} and {\tt Gurobi}
is reported in Table~\ref{tb:compare:61}.
For the Cholesky-Newton method, the Newton direction does not exist
since $s(\lambda) \coloneqq -(H + \lambda I)^{-1}g = \mathbf{0}_n$
for all $\lambda > -\lambda_{\min}$,
and it terminated in the very first loop.
For {\tt fminunc}, we tested it on three initial points:
$s^{(1)} = \mathbf{0}_5$, $s^{(2)} = \mathbf{1}_5$, $s^{(3)} = -10 \cdot \mathbf{1}_5$.
For all of them, {\tt fminunc} converged to
a critical point that is not a global minimizer.
{\tt Baron} found a point whose objective value is
$-5.2479$ in $0.50$ second,
and returned the lower bound $-16.2604$ after $300$ seconds.
{\tt Gurobi} obtained such a point in $0.10$ second,
and returned the lower bound $-5.3516$  after $300$ seconds.
Both {\tt Baron} and {\tt Gurobi} cannot verify the global optimality.
\begin{table}[htbp]
\centering
\caption{Performance comparison of various solvers for Example \ref{ex:6-1}.}
\label{tb:compare:61}
\begin{tabular}{lrlr}
\specialrule{.2em}{0em}{0.1em}
Method & \multicolumn{2}{c}{$M_{value}$} & Time (s) \\
\specialrule{.1em}{.1em}{0.1em}
SDP relax. & $-5.2479$ &   & $0.04$ \\
\specialrule{.1em}{.1em}{0.1em}
Cho.-New. & \multicolumn{2}{c}{not convergent} \\
\specialrule{.1em}{.1em}{0.1em}
\multirow{3}{*}{{\tt fminunc}}
 &  0.0000 & ($s^{(1)}$)   & 0.15 \\
 &  0.0000 & ($s^{(2)}$)   & 0.11 \\
 &  $-0.4999$ & ($s^{(3)}$)  & 0.03 \\
\specialrule{.1em}{.1em}{0.1em}
& Upper Bound & Lower Bound \\
{\tt Baron}  & $-5.2479$ & $-16.2604$ & $\ge 300.00$ \\
{\tt Gurobi}  & $-5.2479$ & $-5.3516$ & $\ge 300.00$ \\
\specialrule{.2em}{0em}{0.1em}
\end{tabular}
\end{table}
\end{Example}

\begin{Example}  \label{ex:6-2}\rm
Consider the CQR problem
\[
\min_{s\in \re^3} \quad   \ve s\ve^4-10\ve s\ve^3
-(\frac{5}{2}s_1^2+2s_2^2+3s_3^2)+s_1s_2+2s_2s_3 + \sum_{i=1}^{3}is_i.
\]
By solving  \reff{problem:dual:sos} and \reff{pro-moment},
we get $\rank\, Y^* =\rank\, Z_1^* = \rank\, Z_2^* = 1$.
It took around $0.01$ second.
By Proposition~\ref{prop:rankone}, we get the minimizer
\[
s^* \, = \,  \big( (Y^*)_{10}, (Y^*)_{20}, (Y^*)_{30} \big) \, \approx \,
(-1.8131,\,3.6458,\, -6.5873).
\]
By Algorithm~\ref{alg:finds}, we get $z^* \approx  7.7441$
and the matrix $R\in \re^{3\times 4}$ has rank $3$.
It confirms $s^*$ is the unique minimizer.
The optimal values are \[ \gamma^* = \vartheta^* =  \mu^* \approx  -1281.5926. \]
The relaxations \reff{problem:dual:sos} and \eqref{pro-moment} are tight.
The absolute and relative errors are
\[
 \texttt{err-abs}  =  2.12 \cdot 10^{-6} , \quad
 \texttt{err-rel}  =  1.65 \cdot 10^{-9}.
\]
Performance of the Cholesky-Newton method,
{\tt fminunc}, {\tt Baron} and {\tt Gurobi}
is reported in Table~\ref{tb:compare:62}.
The Cholesky-Newton method took around $0.01$ second to find the global minimizer.
For {\tt fminunc}, we tested it on three initial points:
$s^{(1)} = \mathbf{0}_3$, $s^{(2)} = \mathbf{1}_3$, $s^{(3)} = (-0.5,0,0.5)$.
It found the global minimizer with the initial point $s^{(1)}$,
while it converged to a critical point that is not
a global minimizer for $s^{(2)}$ and $s^{(3)}$.
{\tt Baron} took around $0.25$ second to
find a feasible point whose objective value is $-1281.5926$,
and returned the lower bound $-1446.0947$ after $300$ seconds.
{\tt Gurobi} took around $0.84$ second to
find a feasible point whose objective value is $-1281.5926$
and to get the lower bound $-1281.6874$.
Both {\tt Baron} and {\tt Gurobi} cannot verify the global optimality.
\begin{table}[htbp]
\centering
\caption{Performance comparison of various solvers for Example~\ref{ex:6-2}.}
\label{tb:compare:62}
\begin{tabular}{lrlr}
\specialrule{.2em}{0em}{0.1em}
Method & \multicolumn{2}{c}{$M_{value}$} & Time (s) \\
\specialrule{.1em}{.1em}{0.1em}
SDP relax. & {$-1281.5926$} & & $0.01$ \\
\specialrule{.1em}{.1em}{0.1em}
Cho.-New. & {$-1281.5926$} & & $0.01$ \\
\specialrule{.1em}{.1em}{0.1em}
\multirow{4}{*}{{\tt fminunc}}
 &  {$-1281.5926$} & ($s^{(1)}$)   & 0.25 \\
 &  {$-1257.7205$} & ($s^{(2)}$)   & 0.06 \\
 &  {$-1257.7205$} & ($s^{(3)}$)  & 0.09 \\
 \specialrule{.1em}{.1em}{0.1em}
& Upper Bound & Lower Bound \\
{\tt Baron}  & $-1281.5926$ & $-1446.0947$ & $\ge 300.00$ \\
{\tt Gurobi}  & $-1281.5926$ & $-1281.6874$ & $0.84$ \\
\specialrule{.2em}{0em}{0.1em}
\end{tabular}
\end{table}
\end{Example}

\begin{Example} \label{ex:6-3} \rm
Consider the CQR problem
\[
\begin{array}{cl}
\min\limits_{s\in \re^5}  & \ve s\ve^4-5\ve s\ve^3-2s_1^2-s_3^2-\frac{3}{2}s_4^2-\frac{1}{2}s_5^2
-2s_1s_2-s_1s_3 \\
& \qquad -3s_1s_4+s_2s_3-2s_2s_4-s_2s_5-2s_3s_5-s_4s_5
+2\sum_{i=1}^{5}s_i.
\end{array}
\]
By solving  \reff{problem:dual:sos} and \reff{pro-moment},
we get $\rank\, Y^* =\rank\, Z_1^* = \rank\, Z_2^* = 1$.
It took around $0.03$ second.
By Proposition~\ref{prop:rankone}, we get the minimizer
\[
s^* \, = \, \big( (Y^*)_{10}, \ldots, (Y^*)_{50} \big) \, \approx \,
(-2.8277,\,-1.4802,\,-0.7917,\, -2.5252,\,-0.9839).
\]
By Algorithm~\ref{alg:finds}, we get $z^*\approx 4.2612 $
and the matrix $R\in \re^{5\times 6}$ has rank $5$.
This confirms $s^*$ is the unique minimizer.
The optimal values are
\[ \gamma^* = \vartheta^* =  \mu^* \approx  -144.8805. \]
The relaxations \reff{problem:dual:sos} and \eqref{pro-moment} are tight.
The absolute and relative errors are
\[
\texttt{err-abs}  = 3.25 \cdot 10^{-7} , \quad
\texttt{err-rel}  =  2.24 \cdot 10^{-9}.
\]
Performance of the Cholesky-Newton method, {\tt fminunc}, {\tt Baron} and {\tt Gurobi}
is reported in Table~\ref{tb:compare:64}.
The Cholesky-Newton method took around $0.01$ second to find the global minimizer.
For {\tt fminunc}, we tested it on three initial points:
$s^{(1)} = \mathbf{0}_5$, $s^{(2)} = \mathbf{1}_5$, $s^{(3)} = (1,0,0,0,0)$.
It found the global minimizer with the initial point $s^{(1)}$,
while it converged to a critical point that is not
a global minimizer for other initial points.
{\tt Baron} took around $0.39$ second to
find a feasible point whose objective value is $-144.8805$,
and returned the lower bound $-2714.0129$ after $300$ seconds.
{\tt Gurobi} took around $1.32$ seconds to
find a feasible point whose objective value is $-144.8805$
and to get the lower bound $-144.8923$.
Both {\tt Baron} and {\tt Gurobi} cannot verify the global optimality.
\begin{table}[htbp]
\centering
\caption{Performance comparison of various solvers for Example \ref{ex:6-3}.}
\label{tb:compare:64}
\begin{tabular}{lrlr}
\specialrule{.2em}{0em}{0.1em}
Method & \multicolumn{2}{c}{$M_{value}$} & Time (s) \\
\specialrule{.1em}{.1em}{0.1em}
SDP relax. & $-144.8805$ & & $0.01$ \\
\specialrule{.1em}{.1em}{0.1em}
Cho.-New. & $-144.8805$ & & $0.01$ \\
\specialrule{.1em}{.1em}{0.1em}
\multirow{3}{*}{{\tt fminunc}}
 &  $-144.8805$ & ($s^{(1)}$)   & 0.15 \\
 &  $-112.6193$ & ($s^{(2)}$)   & 0.11 \\
 &  $-112.6193$ & ($s^{(3)}$)   & 0.12 \\
 \specialrule{.1em}{.1em}{0.1em}
& Upper Bound & Lower Bound \\
{\tt Baron}  & $-144.8805$ & $-2714.0129$ & $\ge 300.00$ \\
{\tt Gurobi}  & $-144.8805$ & $-144.8923$ & $1.32$ \\
\specialrule{.2em}{0em}{0.1em}
\end{tabular}
\end{table}
\end{Example}

\begin{Example}\label{ex:sec6:add:ug=0}\rm
Consider the CQR problem
\[
\min_{s\in \re^3} \quad   \ve s\ve^4-
10s_1^2 + 2s_2^2 + 4s_3^2 - 12\sqrt{2}s_2 - 28\sqrt{2}s_3.
\]
By solving \reff{problem:dual:sos} and \reff{pro-moment},
we get $\rank\, Y^* = 2$, $\rank\, Z_1^* = \rank\, Z_2^*=1$.
It took around 0.03 second. By Algorithm~\ref{alg:finds},
we get $\ve s^*\ve =2.2361$ and two global minimizers
\[
s^{*} \approx (\pm 1.5812,\, 0.7071,\, 1.4142).
\]
The optimal values are \[ \gamma^* = \vartheta^* =  \mu^* \approx  -59.0000. \]
The absolute and relative errors are
\[
\texttt{err-abs}  = 2.95 \cdot 10^{-7} , \quad
\texttt{err-rel}  =  5.00 \cdot 10^{-9}.
\]
Performance of the Cholesky-Newton method, {\tt fminunc}, {\tt Baron} and {\tt Gurobi}
is reported in Table~\ref{tb:compare:sec6:add:ug=0}.
Note that the eigenvector of $H$ associated to
$\lmd_{min}$ is $\bbm 1 & 0  & 0 \ebm^T$, which is orthogonal to the vector $g$,
and the matrix $B(s^*)$ is singular.
The Cholesky-Newton method cannot find the global minimizer.
For {\tt fminunc}, we tested it on three initial points:
$s^{(1)} = \mathbf{0}_3$, $s^{(2)} = \mathbf{1}_3$, $s^{(3)} = 10\cdot \mathbf{1}_3.$
It found the global minimizer with the initial points $s^{(2)}$ and $s^{(3)}$,
while it converged to a critical point for the initial point $s^{(1)}$.
{\tt Baron} took around $0.41$ seconds to find the global minimizer.
{\tt Gurobi} took around $0.16$ seconds to
find a feasible point whose objective value is $-59.0000$
and to get the lower bound $-59.0039$.
It cannot verify the global optimality.
\begin{table}[htbp]
\centering
\caption{Performance comparison of various solvers for Example \ref{ex:sec6:add:ug=0}.}
\label{tb:compare:sec6:add:ug=0}
\begin{tabular}{lrlr}
\specialrule{.2em}{0em}{0.1em}
Method & \multicolumn{2}{c}{$M_{value}$} & Time (s) \\
\specialrule{.1em}{.1em}{0.1em}
SDP relax. & $-59.0000$ &  & $0.03$ \\
\specialrule{.1em}{.1em}{0.1em}
Cho.-New. & \multicolumn{2}{c}{not convergent} & \\
\specialrule{.1em}{.1em}{0.1em}
\multirow{3}{*}{{\tt fminunc}}
 &  $-55.7283$ & ($s^{(1)}$)   & 0.06 \\
 &  $-59.0000$ & ($s^{(2)}$)   & 0.04 \\
 &  $-59.0000$ & ($s^{(3)}$)  & 0.02 \\
\specialrule{.1em}{.1em}{0.1em}
& Upper Bound & Lower Bound \\
{\tt Baron}  & $-59.0000$ & $-59.0000$ & $0.41$  \\
{\tt Gurobi}  & $-59.0000$ & $-59.0039$ & $0.16$ \\
\specialrule{.2em}{0em}{0.1em}
\end{tabular}
\end{table}
\end{Example}

\begin{Example}\rm
\label{ex:sec6:add:ug=0:pos}
Consider the CQR problem
\[
\min_{s\in \re^3} \quad   \ve s\ve^4-
\frac{8}{3}\|s\|^3+ 2s_2^2 + 4s_3^2-4s_2.
\]
By solving \reff{problem:dual:sos} and \reff{pro-moment}, we get
$\rank\, Y^* = 2$, $\rank\, Z_1^* = \rank\, Z_2^*=1$.
It took around $0.01$ second. By Algorithm~\ref{alg:finds}, we get $\ve s^* \ve =2$,
\[
 s^{*} \, = \, (\pm \sqrt{3},\, 1,\, 0).
\]
The optimal values are $\gamma^* = \vartheta^* =  \mu^* \approx  -7.3333$.
The absolute and relative errors are
\[
\texttt{err-abs}  = 1.70\cdot 10^{-7} , \quad
\texttt{err-rel}  = 2.32 \cdot 10^{-8}.
\]
Performance of the Cholesky-Newton method, {\tt fminunc}, {\tt Baron} and {\tt Gurobi}
is reported in Table~\ref{tb:compare:sec6:add:ug=0:pos}.
For this CQR, the eigenvector of $H$ associated to
$\lmd_{min}$ is $\bbm 1 & 0  & 0 \ebm^T$,
which is orthogonal to the vector $g$, and the matrix $B(s^*)$ is singular.
The Cholesky-Newton method cannot find the global minimizer.
For {\tt fminunc}, we tested it on three initial points:
$s^{(1)} = \mathbf{0}_3$, $s^{(2)} = \mathbf{1}_3$,
$s^{(3)} = 10 \cdot \mathbf{1}_3$.
It found the global minimizer with the initial points $s^{(2)}$ and $s^{(3)}$,
while it converged to a critical point that is not
a global minimizer for the initial point $s^{(1)}$.
{\tt Baron} took around $0.16$ second to
find a feasible point whose objective value is $-7.3333$,
and returned the lower bound $-7.6446$ after $300$ seconds.
It cannot verify the global optimality.
{\tt Gurobi} took around $0.31$ second to find the global minimizer.
\begin{table}[htbp]
\centering
\caption{Performance comparison of various solvers for
Example~\ref{ex:sec6:add:ug=0:pos}.}
\label{tb:compare:sec6:add:ug=0:pos}
\begin{tabular}{lrlr}
\specialrule{.2em}{0em}{0.1em}
Method & \multicolumn{2}{c}{$M_{value}$} & Time (s) \\
\specialrule{.1em}{.1em}{0.1em}
SDP relax. & $-7.3333$ &  & $0.01$ \\
\specialrule{.1em}{.1em}{0.1em}
Cho.-New. & \multicolumn{2}{c}{not convergent} & \\
\specialrule{.1em}{.1em}{0.1em}
\multirow{3}{*}{{\tt fminunc}}
 &  $-5.7879$ & ($s^{(1)}$)   & 0.03 \\
 &  $-7.3333$ & ($s^{(2)}$)   & 0.00 \\
 &  $-7.3333$ & ($s^{(3)}$)  & 0.01 \\
\specialrule{.1em}{.1em}{0.1em}
& Upper Bound & Lower Bound \\
{\tt Baron}  & $-7.3333$ & $-7.6446$ & $\ge 300.00$ \\
{\tt Gurobi}  & $-7.3333$ & $-7.3333$ & $0.31$ \\
\specialrule{.2em}{0em}{0.1em}
\end{tabular}
\end{table}
\end{Example}

In the following, we explore the computational performance of
the SDP relaxations \reff{problem:dual:sos} and \reff{pro-moment}.
Up to a scaling to the variable $s=  (\frac{4}{\sigma})^{1/4} \tilde{s}$,
the CQR problem can be transformed to a new one with $\sigma = 4$,
i.e., the coefficient of $\ve s \ve^4$ is one.

\begin{Example}\label{ex:random} \rm
Consider some randomly generated CQR problem in the form
\be \label{CQR:rand}
\min_{s\in \re^n} \quad  \ve s\ve^4 + \frac{\bt}{6} \ve s\ve ^3
+\frac{1}{2}s^THs+g^Ts,
\ee
where $\sigma = 4$. The vector $g$ and matrix $H$ are generated
in {\tt MATLAB} as in \cite{zhu2025CQR}:
\[
g = \texttt{randn(n,1)}, \quad  H_1  = \texttt{randn(n,n)},
\quad H = \frac{1}{2}\big( H_1 + H_1^T \big).
\]
For each $n$ and $\beta$, we generate $20$ instances. For each instance,
we solve the CQR by SDP relaxations \reff{problem:dual:sos} and \reff{pro-moment}.
For all instances, the CQR problems are solved successfully
by \reff{problem:dual:sos} and \reff{pro-moment}.
The accuracy error $\texttt{err-abs}$ is around $10^{-9}$ for all instances.
We report the average computational time (in seconds) for each case of $(n, \beta)$.
The numerical performance is reported in Table~\ref{tb:random}.
As we can see, these CQR problems can be solved efficiently by the relaxations.
For instance, when the dimension $n=1000$,
the relaxations can be solved within around one minute.
\begin{table}[htbp]
\centering
\caption{Computational time (in seconds) for solving \reff{CQR:rand}
by \reff{problem:dual:sos} and \reff{pro-moment}. }
\label{tb:random}
\begin{tabular}{c|rrrrrr}
\specialrule{.2em}{0em}{0.1em}
\diagbox{$n$}{$\beta$} & 10 & 1 & 0 & -1 & -10 & -100\\
\specialrule{.1em}{.1em}{0.1em}
100 & 0.11 & 0.12 & 0.11 & 0.12 & 0.11 & 0.24\\
\specialrule{.1em}{.1em}{0.1em}
200 & 0.34 &  0.30 & 0.30 & 0.30 & 0.45 & 0.97 \\
\specialrule{.1em}{.1em}{0.1em}
300 & 1.72 & 1.76 & 1.74 & 1.77 & 1.63 & 3.49\\
\specialrule{.1em}{.1em}{0.1em}
400 & 5.01 & 5.84 & 5.79 & 5.96 & 5.20 & 11.30\\
\specialrule{.1em}{.1em}{0.1em}
500 & 6.15 & 6.19 & 6.27 & 6.21 & 6.45 & 21.82\\
\specialrule{.1em}{.1em}{0.1em}
600 & 11.60 & 11.97 & 11.54 & 11.30 & 11.30 & 35.78\\
\specialrule{.1em}{.1em}{0.1em}
700 & 18.72 & 18.43 & 18.07 & 18.82 & 18.84 & 43.14\\
\specialrule{.1em}{.1em}{0.1em}
800 & 23.64 & 23.98 & 22.91 & 23.85 & 23.36 & 79.03\\
\specialrule{.1em}{.1em}{0.1em}
900 & 37.09 & 35.87 & 36.54 & 34.78 & 35.54 & 67.53 \\
\specialrule{.1em}{.1em}{0.1em}
1000 & 52.92 & 52.74 & 53.46 & 54.15 & 54.33 & 60.93\\
\specialrule{.2em}{0em}{0.1em}
\end{tabular}
\end{table}
Moreover, we also test the Cholesky-Newton method and {\tt Gurobi} for these instances. By the way of generation, the matrix $H$ is usually not positive definite,
and $g^Tv$ is generally nonzero for the eigenvector $v$ of $H$
corresponding to $\lmd_{min}$. So, this CQR problem generally belongs
to Case~1 or 2 of \cite[Theorem~2.5]{zhu2025CQR}.
For all instances, the Cholesky-Newton method found global minimizers for CQR subproblems, and it took less than one second
to find the global minimizer when $n=1000$ and $\beta = -100$,
which demonstrates good scalability and computational efficiency
compared with the SDP relaxation method
for regular CQR problems such as (\ref{CQR:rand}).
For {\tt Gurobi}, it does not converge within $300$ seconds
for all instances with $n=100$ and $\beta = -1$.
{\tt Baron} is not tested, since the problem size exceeds
the limit of its publicly available version.
\end{Example}

In Table~\ref{tab:cqr-example-summary},
we summarize previous CQR examples,
including their defining parameters, global minimizers, and tightness of SDP relaxations.
The case $g^Tv=0$, where $v$ is an eigenvector of $H$ associated with the minimum eigenvalue $\lambda_{\min}$, belongs to the ``hard case" for the secular equation approach.
CQR problems can be classified into five cases
as in \cite[Theorem~2.5]{zhu2025CQR}.

\renewcommand{\arraystretch}{1.35} 
\begin{table}[htbp]  
\centering
\scriptsize 
\caption{Summary of some CQR examples.}
\label{tab:cqr-example-summary}
\begin{tabular}{@{}c c c p{0.20\textwidth} p{0.22\textwidth} p{0.20\textwidth}@{}}
\toprule
Exam. & $n$ & $(\beta,\sigma)$ & $H$ and $g$
& CQR classification 
& SDP tightness \\
\midrule

\ref{ex:cases} (i)
& $2$
& $(-16,4)$
& $\displaystyle H=\operatorname{diag}(4,2)$,
  \newline
  $\displaystyle g=(-8,0)$
& Case~4, unique minimizer
  $s^*=(2,0)$.
& (\ref{bt>=-3}) holds; tight.
\newline
  $\beta+3\sigma\|s^*\|=8$. \\[1.0em]
  \hline

\ref{ex:cases} (ii)
& $2$
& $(-36,4)$
& $\displaystyle H=26I_2$,
  \newline
  $\displaystyle g=(12,0)$
& Case~3, two minimizers
  $(-1,0)$ and $(-2,0)$, with different norms.
& (\ref{bt>=-3}) fails; not tight.
\newline
  $\beta+3\sigma\|s^*\|=-24$ or
  $-12$,
\newline
  $\vartheta^*=-5<\mu^*=0$. \\[1.0em]
  \hline

\ref{ex:cases} (iii)
& $2$
& $(-24,4)$
& $\displaystyle H=\operatorname{diag}(20,10)$,
  \newline
  $\displaystyle g=(-12,-2)$
& Case~5, global minimum value obtained by Moment-SOS.
& (\ref{bt>=-3}) fails; not tight.
\newline
 $\beta+3\sigma\|s^*\|\approx -2.36$;
\newline
$\vartheta^*=-6<\mu^*\approx -5.64$. \\[1.0em]
\hline

\ref{ex:cases} (iv)
& any
& \makecell[t]{
$(-12+6\epsilon,4)$, \\
$\epsilon\in(0,2)$
}
& $\displaystyle H=2I_n$,
  \newline
  $\displaystyle g={\bf0}_n$ (hard case)
& Case~3 or 5, depending on $\epsilon$;
  unique minimizer $s^*={\bf0}_n$.
& (\ref{bt>=-3}) holds; tight.
\newline
$\|s^*\|=0$. \\[1.0em]
\hline

\ref{ex:cases} (v)
& any
& $(-24,4)$
& $\displaystyle H=8I_n$,
  \newline
  $\displaystyle g={\bf0}_n$ (hard case)
& Case~4, minimizers are
  ${\bf0}_n$ and all $s^*$ with
  $\|s^*\|=2$;
  $B(s^*)$ is singular for all nonzero $s^*$.
& (\ref{bt>=-3}) holds; tight.
\newline
$\beta+3\sigma\|s^*\|=0$ for all nonzero $s^*$. \\[1.0em]
\hline

\ref{ex:6-1}
& $5$
& $(-6,4)$
& $\displaystyle H=\mathbf{1}_5\mathbf{1}_5^{T}-6I_5$,
\newline
$\displaystyle g={\bf0}_n$ (hard case)
& Case~1, infinitely many minimizers with
$\|s^*\|\approx1.6559$ and singular $B(s^*)$.
& (\ref{bt>=-3}) holds; tight.
\newline
$\beta+3\sigma\|s^*\|\approx13.87$. \\[1.0em]
\hline

\ref{ex:6-2}
& $3$
& $(-60,4)$
& $\displaystyle H\prec0$,
\newline
$\displaystyle g=(1,2,3)$
& Case~1, unique minimizer.
& (\ref{bt>=-3}) holds; tight.
\newline
$\beta+3\sigma\|s^*\|\approx32.93$. \\[1.0em]
\hline

\ref{ex:6-3}
& $5$
& $(-30,4)$
& $H$ is indefinite,
\newline
$\displaystyle g=(2,2,2,2,2)$
& Case~1, unique minimizer.
& (\ref{bt>=-3}) holds; tight.
\newline
$\beta+3\sigma\|s^*\|\approx21.13$. \\[1.0em]
\hline
\ref{ex:sec6:add:ug=0}
& $3$
& $(0,4)$
& $H$ is indefinite,
\newline
$\displaystyle g=-\sqrt2(0,12,28)$
\newline
(hard case since $g^Tv=0$).
& Case~1, two minimizers with
$\|s^*\|\approx2.2361$ and singular $B(s^*)$.
& (\ref{bt>=-3}) holds; tight.
\newline
$\beta+3\sigma\|s^*\|\approx26.83$. \\[1.0em]
\hline

\ref{ex:sec6:add:ug=0:pos}
& $3$
& $(-16,4)$
& $\displaystyle H=\operatorname{diag}(0,4,8)$,
\newline
$\displaystyle g=(0,-4,0)$
\newline
(hard case since $g^Tv=0$)
& 
two minimizers $s^*=(\pm\sqrt3,1,0)$
with singular $B(s^*)$.
& (\ref{bt>=-3}) holds; tight.
\newline
$\beta+3\sigma\|s^*\|=8$. \\
\hline

\ref{ex:random}
& \makecell[t]{
$100$\\
$\sim$\\
$1000$}
&
\makecell[t]{
$\beta=-100\sim10$\\
$\sigma=4$
}
& $H=\frac{1}{2}(H_1+H_1^T)$
\newline
$H_1=\texttt{randn(n,n)}$
\newline
$g=\texttt{randn(n,1)}$
& $H$ is typically not pd.
& (\ref{bt>=-3}) holds; tight. \\

\bottomrule
\end{tabular}   
\end{table}

\subsection{Usage of SDP relaxations (\ref{problem:dual:sos}) and (\ref{pro-moment}) in the CQR scheme}

In the following, we explore how to apply SDP relaxations (\ref{problem:dual:sos}) and (\ref{pro-moment}) to solve CQR subproblems that appear in numerical methods
for solving the unconstrained optimization problem (\ref{minf(x)}).
For convenience of description, we follow \cite[Algorithm~2]{zhu2025CQR}
to formulate CQR subproblems and then solve them by (\ref{problem:dual:sos}) and (\ref{pro-moment}).
We let $x^{(k)}$ denote the point in the $k$th outer iteration,
and follow the parameter settings and the termination criterion given as in \cite[Section~4.1.1]{zhu2025CQR}.
The parameter $\beta$ is given by Option~1 of \cite[Algorithm~5]{zhu2025CQR}.
In the experiments, all CQR subproblems are solved successfully by
the SDP relaxations (\ref{problem:dual:sos}) and (\ref{pro-moment}).

\begin{Example}\rm \label{ex:ext_Rosenbrock}
Consider the Extended Rosenbrock function \cite{More1981testing} in $x \in \re^{100}$:
\begin{equation}  \label{fun:extRos}
f(x)  \, = \,  \sum_{i=1}^{50} \left[ 100(x_{2i} - x_{2i-1}^2)^2
+ (1 - x_{2i-1})^2 \right],
\end{equation}
Its minimum value is $0$ and the global minimizer is $x^* = {\bf 1}_{100}$.
We implement the CQR scheme \cite[Algorithm~2]{zhu2025CQR}
to minimize $f(x)$, with the following initial points:
\[
\begin{gathered}
u^{(1)} = {\bf 0}_{100},\quad u^{(2)} = -{\bf 1}_{100},\quad
u^{(3)} = (-1.2,\, 1.0, \,\dots, -1.2,\, 1.0),	
\end{gathered}
\]
where $u^{(3)}$ is the initial point suggested in \cite{More1981testing}.
The SDP relaxations (\ref{problem:dual:sos}) and (\ref{pro-moment})
are applied to solve the appearing CQR subproblems.
For all these three initial points, the CQR scheme converged to the global minimizer $x^*$. For the initial point $u^{(1)}$, the termination criterion is achieved after $260$ outer iterations, while it took $272$ outer iterations for the initial point $u^{(2)}$ and $338$ outer iterations for the initial point $u^{(3)}$.
The consumed time is $28.98$ seconds, $27.45$ seconds, and $32.81$ seconds, respectively. The objective function value $f(x^{(k)})$
versus the outer iteration $k$ is plotted in Figure~\ref{fig:rosen100}.
For better visualization, the vertical axis for the objective function value
is shown in the logarithmic scale with base $10$.
We also use the Newton-Cholesky method \cite[Algorithm~1]{zhu2025CQR}
to solve the CQR subproblems for minimizing \reff{fun:extRos}.
For all three initial points, the CQR scheme successfully converged
to the global minimizer $x^*$ of $f(x)$ after $260$, $272$, and $338$ outer iterations,
with time consumptions of $0.43$, $0.41$ and $0.43$ second respectively.
In contrast, {\tt Gurobi} and \texttt{fmincon}
failed to solve the CQR subproblems globally,
since they are highly nonconvex.
{\tt Baron} is not tested,
since the problem size exceeds the limit of its publicly available version.
\begin{figure}[htbp]
\centering
\caption{Convergence of the CQR scheme with the SDP relaxations (\ref{problem:dual:sos}) and (\ref{pro-moment}) on the Extended Rosenbrock function.}
\includegraphics[width=0.67\textwidth]{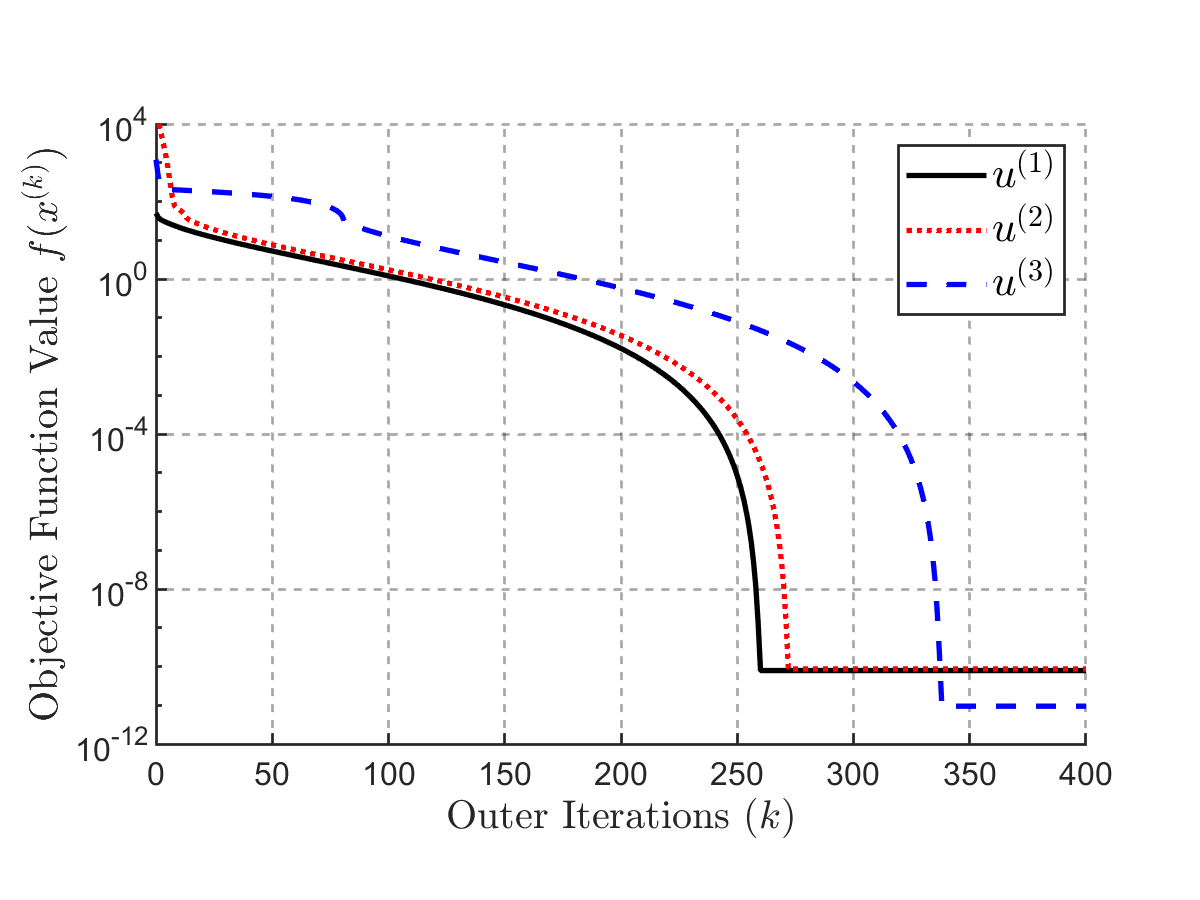}
\label{fig:rosen100}
\end{figure}
\end{Example}

For solving nonconvex nonlinear optimization problems,
a major difficulty is how to escape from critical points that are not global minimizers.
This is particularly challenging for the CQR scheme.
When $x^{(k)}$ is a critical point that is not a global minimizer,
we have $g(x^{(k)}) \coloneqq \nabla f(x^{(k)}) = 0$,
so finding global minimizers for the CQR subproblem from secular equations may be difficult. However, the SDP relaxations (\ref{problem:dual:sos}) and (\ref{pro-moment})
can still solve CQR subproblems globally when $g(x^{(k)}) = 0$.
Therefore, if one applies the relaxations (\ref{problem:dual:sos}) and (\ref{pro-moment})
to globally solve the appearing CQR subproblems, then the adaptive regularization method
can escape from critical points that are not global minimizers of (\ref{minf(x)}).
The following example is such an exposition.

\begin{Example} \label{ex:ar3:com:n=5}\rm
Consider the nonlinear function in $x \in \re^{100}$:
\[
f(x) = \sum_{i=1}^{100} \left( \frac{1}{4} x_i^4 - \frac{1}{2} x_i^2 \right)
+ \frac{1}{2} \sum_{i=1}^{99} (x_i + x_{i+1})^2 + \frac{1}{2}(x_1 + x_{100})^2.
\]
Its minimum value is $-25$ and the global minimizers are
\[
x^* \, = \, \pm (1,-1,1,-1, \ldots, 1,-1).
\]
We implement the CQR scheme \cite[Algorithm~2]{zhu2025CQR}
with three initial points:
\[
u^{(1)} = {\bf 0}_{100},\quad u^{(2)} = {\bf 1}_{100},\quad
u^{(3)} = -\frac{1}{2}{\bf 1}_{100},
\]
where $u^{(1)} = {\bf 0}_{100}$ is a critical point of $f$ that is not a minimizer.
The SDP relaxations (\ref{problem:dual:sos}) and (\ref{pro-moment})
are applied to solve each CQR subproblem. For all three initial points,
the CQR scheme converged to the global minimizer $x^*$.
For the initial point $u^{(1)}$,
the termination criterion is met after $15$ outer iterations,
while it took $17$ outer iterations for
$u^{(2)}$ and $17$ outer iterations for $u^{(3)}$.
For all three initial points, the global minimum $-25$ is obtained.
The time consumptions are $1.01$, $1.17$ and $1.05$ seconds respectively.
The objective function value $f(x^{(k)})$ versus the outer iteration $k$
are plotted in Figure~\ref{fig:saddle100}.
We also use the Newton-Cholesky method \cite[Algorithm~1]{zhu2025CQR}
to solve the CQR subproblems. For the initial point $u^{(1)}$,
the algorithm terminated after $1$ iteration
at the saddle point with the objective value $0$,
failing to find the global minimizer.
For the initial points $u^{(2)}$ and $u^{(3)}$,
the CQR scheme successfully converged to $x^*$
after $17$ outer iterations, with time consumptions of
$0.20$ second and $0.16$ second respectively.
Moreover, \texttt{fmincon} failed to escape the saddle point $u^{(1)}$,
terminating after $1$ iteration with the objective value $0$.
For the initial points $u^{(2)}$ and $u^{(3)}$,
\texttt{fmincon} converged to a saddle point
that is not a global minimizer, after $11$ and $9$ iterations respectively,
with the final objective value $0$.
{\tt Gurobi} failed to find global minimizers of CQR subproblems
and the CQR scheme did not converge to a global minimizer of $f(x)$.
{\tt Baron} is not tested, since
the problem size exceeds the limit of its publicly available version.
\begin{figure}[htbp]
\centering
\caption{Convergence of the CQR scheme with the SDP relaxations (\ref{problem:dual:sos}) and (\ref{pro-moment}) on Example~\ref{ex:ar3:com:n=5}.}
\includegraphics[width=0.67\textwidth]{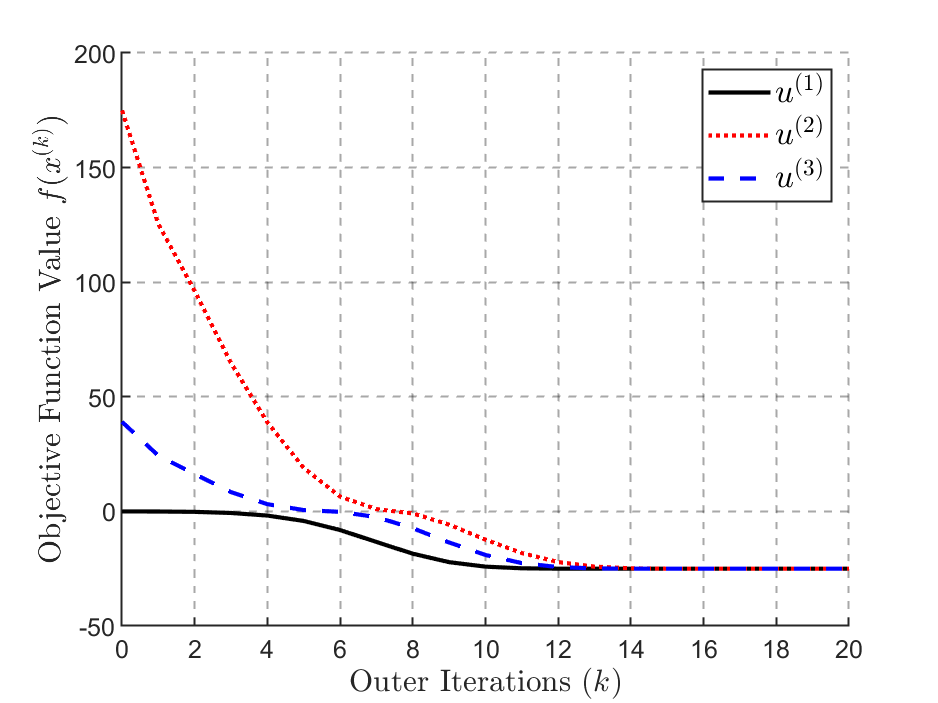}
\label{fig:saddle100}
\end{figure}
\end{Example}

\section{Conclusions}
\label{sc:con}

In this paper, we give an SDP relaxation method for
finding global minimizers of the cubic-quartic regularization problem.
This method transforms the original non-convex optimization problem into
a tractable semidefinite program.
We prove that the relaxations are tight
when the regularization parameters satisfy
\[ \ve s^* \ve (\beta + 3\sigma \ve  s^*\ve) \, \ge \, 0, \]
where $s^*$ is a global minimizer of the CQR problem.
This condition may be used as guidelines for parameter selections
in numerical optimization methods.
Moreover, we show that if the relaxations are tight,
then all nonzero minimizers of the CQR problem have the same length.
Algorithm \ref{alg:finds} is given
to obtain the set of all global minimizers, when the SDP relaxations are tight.
Numerical experiments demonstrate that our SDP relaxation method
is efficient for solving CQR problems.

The SDP relaxations (\ref{problem:dual:sos}) and (\ref{pro-moment})
can be solved in $O(n^{3.5} \ln(1/\eps))$ arithmetic operations
by interior-point methods \cite{Todd01}, so they are polynomial time solvable.
Therefore, for the case that (\ref{problem:dual:sos}) and (\ref{pro-moment})
are tight, the CQR problems are polynomial time solvable.
This is the case when $H$ has a nonpositive eigenvalue or $\beta \ge 0$.
We remark that this is the first method that can solve
a class of CQR problems in polynomial time, to the best of the authors' knowledge.

When SDP relaxations (\ref{problem:dual:sos}) and (\ref{pro-moment}) are not tight,
Algorithm~\ref{alg:finds} cannot produce global minimizers for the CQR problem (\ref{CQR:POLY}).
For these cases, solving the CQR problem is more difficult.
One can apply Moment-SOS relaxation methods to solve
the equivalent polynomial optimization problem \reff{eq:trans_pop}.
See Remark~\ref{rmk:pop_reform} and Example~\ref{ex:cases}(iii).
The problem \reff{eq:trans_pop} can be solved by tight
Moment-SOS relaxations; see the work \cite{NDS06,Tight18}.
Typically, Moment-SOS relaxations are expensive for large scale problems.
It is interesting future work to develop practically efficient methods
for solving large CQR problems when the relaxations
\reff{problem:dual:sos} and \reff{pro-moment} are not tight.

Finally, we remark that the Newton-Cholesky method based on secular equations is more scalable for CQR problems not belonging to the hard cases, as demonstrated in Example~\ref{ex:random}.
On the other hand, the SDP relaxations (\ref{problem:dual:sos}) and (\ref{pro-moment}) can find all global minimizers under tightness conditions, even if the CQR problem belongs to the hard case.
This is particularly useful when the adaptive regularization method needs to escape from a nonglobal critical point; see Example~\ref{ex:ar3:com:n=5}.
Our SDP relaxation method is a certifiable approach for solving CQR problems globally, under tightness conditions.

\bigskip \noindent
{\bf Acknowledgement}
Jiawang Nie is partially supported by the NSF grant DMS-2513254.

\end{document}